\documentclass[12pt,reqno]{amsart}
\usepackage{amsmath,amsfonts,amssymb,amscd,amsthm,amsbsy, color}
\usepackage{graphicx}
\usepackage[usenames]{xcolor}
\usepackage{comment}

\numberwithin{equation}{section}

\newtheorem{meta-thm}[theorem]{Meta-Theorem}

\newcommand\beq[1]{ \begin{equation}\label{#1} }
\newcommand{\eeq}{ \end{equation} }

\newcommand\beqa[1]{ \begin{eqnarray} \label{#1}}
\newcommand{\eeqa}{ \end{eqnarray} }
\newcommand{\beqano}{ \begin{eqnarray*} }
\newcommand{\eeqano}{ \end{eqnarray*} }
\newcommand\equ[1]{{\rm (\ref{#1})}}

\def\E{{\mathcal E}}

\def\integer{{\mathbb Z}}

\def\E{{\mathcal E}}

\def\R{{\mathcal R}}

\begin{document}

\title[Qualitative and analytical results of bifurcations to halo orbits]
{Qualitative and analytical results of the bifurcation thresholds to halo orbits}

\author[S. Bucciarelli]{Sara Bucciarelli}
\address{
Departments of Mathematics, University of Roma Tor Vergata, Via della Ricerca Scientifica 1,
00133 Roma (Italy)}
\email{kayleigh.s@hotmail.it}

\author[M. Ceccaroni]{Marta Ceccaroni}
\address{
Departments of Mathematics, University of Roma Tor Vergata, Via della Ricerca Scientifica 1,
00133 Roma (Italy)}
\email{ceccaron@mat.uniroma2.it}

\author[A. Celletti]{Alessandra Celletti}
\address{
Departments of Mathematics, University of Roma Tor Vergata, Via della Ricerca Scientifica 1,
00133 Roma (Italy)}
\email{celletti@mat.uniroma2.it}

\author[G. Pucacco]{Giuseppe Pucacco}
\address{
Departments of Physics, University of Roma Tor Vergata, Via della Ricerca Scientifica 1,
00133 Roma (Italy)}
\email{pucacco@roma2.infn.it}



\baselineskip=18pt              




\begin{abstract}
We study the dynamics in the neighborhood of the collinear Lagrangian points in the spatial, circular, restricted three--body problem.
We consider the case in which one of the primaries is a radiating body and the other is oblate
(although the latter is a minor effect).
Beside having an intrinsic mathematical interest, this model is particularly suited for the description of a mission of a
spacecraft (e.g., a solar sail) to an asteroid.

The aim of our study is to investigate the occurrence of bifurcations to halo orbits, which
take place as the energy level is varied. The estimate of the bifurcation thresholds is performed by
analytical and numerical methods: we find a remarkable agreement between the two approaches.
As a side result, we also evaluate the influence of the different parameters,
most notably the solar radiation pressure coefficient, on the dynamical behavior of the
model.

To perform the analytical and numerical computations, we start by implementing
a center manifold reduction. Next, we estimate the bifurcation values using qualitative techniques
(e.g. Poincar\'e surfaces, frequency analysis, FLIs).
Concerning the analytical approach, following \cite{CPS} we implement a resonant normal form, we transform to suitable action-angle
variables and we introduce a detuning parameter measuring the displacement from the synchronous resonance.
The bifurcation thresholds are then determined as series expansions in the detuning. Three concrete examples are considered
and we find in all cases a very good agreement between the analytical and numerical results.
\end{abstract}

\subjclass[2010]{37N05, 70F15, 37L10}
\keywords{Collinear Lagrangian points, Center manifolds, Halo orbits}

\maketitle

\tableofcontents

\section{Introduction}\label{intro}
We consider the motion of a small body with negligible mass in the gravitational field of two primaries which move on circular
trajectories around their common barycenter. We refer to this model as the
\sl spatial, circular, restricted \rm three--body problem (hereafter SCR3BP).
As it is well known (see, e.g., \cite{Szebehely}), the SCR3BP admits five equilibrium positions in the synodic reference
frame, which rotate with the angular velocity of the primaries. Two of such positions make an equilateral triangle
with the primaries, while the other three equilibria are collinear with the primaries.

While the triangular positions are shown to be stable for a wide range of the mass ratio of the primaries,
the collinear points are unstable. Nevertheless, the collinear equilibria turn out to be very useful in low--energy space missions
and particular attention has been given to the so--called \sl halo orbits, \rm which are periodic trajectories around the collinear points, generated when increasing the energy level as bifurcations from the so--called \sl planar Lyapunov family \rm of periodic orbits.

We will consider a model in which one of the primaries (e.g. the Sun) is radiating; we will see that in some cases the
effect of the solar radiation pressure (hereafter SRP) on a \sl solar sail, \rm namely an object with a high value of the area-to-mass ratio, is to lower the bifurcation threshold, enabling bifurcations before unfeasible.
For sake of generality, we also consider the case in which the other primary is oblate, although this effect is
definitely negligible in many concrete applications (\cite{SS1}). We will consider the following three case studies.
The first one describes the interaction between the Earth--Moon barycenter and the Sun; we will refer to this case
as the \sl Sun--barycenter \rm system. The second sample is provided by the Earth--Moon system. The last case describes
the interaction between the Sun and one of the largest asteroids, Vesta, for which the effect of SRP is very important.


Our model depends on three main parameters, which are the
mass ratio $\mu$ of the primaries, the performance $\beta$ of the sail providing the ratio between the
acceleration of the radiation pressure to the gravitational acceleration of the main primary (equivalently, one can use the parameter
$q=1-\beta$) and the oblateness denoted by $A$, defined as $A=J_{2}r_{e}^{2}$, see \cite{DM},
where $J_{2}$ is the so--called dynamical oblateness coefficient and $r_e$ is the equatorial radius of the planet.
We remark that the orientation of the sail is set to be perpendicular with respect to the direction joining it with the main primary,
so that the system still preserves the Hamiltonian character and no dissipation is allowed.
All results could be easily generalised to the case in which the orientation
of the sail is kept constant (i.e., non necessarily perpendicular)  with respect to the
direction joining it with the main primary as the system would still keep its Hamiltonian nature.

The aim of this paper is to make concrete analytical estimates of the bifurcations to
halo orbits based on an extension of the theory developed in \cite{CPS} and
to compare the mathematical results with a qualitative investigation of the dynamics.
As a side result, we shall evaluate the role of the parameters of the model, in particular the solar radiation pressure coefficient and
the oblateness.

The collinear points with SRP and oblateness are shown to be of $saddle\times center\times center$ type for typical
parameter values (compare with Sections~\ref{sec:3BP} and \ref{sec:linearstab}). According
to widespread techniques, the dynamics can be conveniently described after
having applied a reduction to the center manifold, which allows us to remove
the hyperbolic components. The procedure consists in applying suitable changes of variables
to reduce the Hamiltonian function to a simpler form and then to compute a Lie series transformation to get rid of the
hyperbolic direction (see, e.g., \cite{JM}, \cite{Ari}). This procedure is performed in Sections~\ref{sec:linearstab} and \ref{sec:center},
taking care of the modifications required by the
consideration of the SRP and of the oblateness. Once the center manifold Hamiltonian has been obtained, we proceed
to implement some numerical methods to investigate the dynamics as the energy level is varied, precisely
we compute some Poincar\'e surfaces of section, we perform a frequency analysis and we determine the Fast Lyapunov Indicators
(\cite{Laskar,LFC,froes}).
The independent and complementary application of the three methods has several advantages, as it helps to unveil
many details of the complicated structures arising around the bifurcation values.

As for the analytical estimates (see \cite{CPS}), after the center manifold reduction we need also to construct
a fourth order normal form around the synchronous resonance.
After introducing a coordinate change to action--angle variables for the quadratic part
of the Hamiltonian, we recognize the existence of a first integral of motion. Finally,
we introduce a quantity called \sl detuning, \rm which measures the (small) discrepancy of the
frequencies around the synchronous resonance. The bifurcation threshold will then be computed
at the first order in the powers series expansion in the detuning (see Section~\ref{sec:analytical}).
We stress that, although we consider just the first non-trivial order of the expansion in the detuning,
the analytical estimates are already in remarkable agreement with the numerical expectation
(see Table~\ref{tab:bif}). Clearly, more refined analytical results can be obtained computing
higher order normal forms, at the expense of a bigger computational effort and provided the results are
performed within the \sl optimal \rm order of normalization.\\

To conclude, let us mention that our comparison is made on the three case studies mentioned before (Sun--barycenter, Earth--Moon, Sun--Vesta); however,
we provide full details only for the Sun--Vesta case, since the other two samples have been already
investigated in the literature with an extensive use of the Poincar\'e maps (see, e.g., \cite{GM,JM}).
The reason to focus on the case of Vesta is the following:
due to the fact that Vesta has a relatively small mass, the interplay of the parameters gives rise to an interesting
and non--trivial dynamical behavior (see Section~\ref{sec:applications}). Indeed, in the Sun--Vesta case we shall see that
the SRP plays an important effect, enabling more bifurcations at relatively low energy levels.

In fact, in the purely gravitational model only the first bifurcation (the \sl standard \rm halo) occurs in general (\cite{GM}), when the  frequency
of the planar Lyapunov orbit is the same as the frequency of its vertical perturbation. A second bifurcation, in which the planar Lyapunov
trajectory regains stability and a second unstable family appears, is a rather extreme phenomenon at high energy values. In the presence of SRP this second bifurcation occurs instead at a much lower value of the energy. Moreover, also a third bifurcation may occur in which the \sl vertical \rm Lyapunov loses stability and the second family disappears. We will show how it is possible to predict the influence of SRP on the thresholds with a first-order perturbation approach.

\vskip.1in

This paper is organized as follows. In Section~\ref{sec:3BP} we introduce the equations of motion and we
determine the equilibrium points, whose linear stability as the parameters are varied is discussed in
Section~\ref{sec:linearstab}, where we prepare the Hamiltonian for the center manifold reduction of
Section~\ref{sec:center}.
The analytical method to determine the bifurcation thresholds is presented in Section~\ref{sec:analytical}, while
in Section~\ref{sec:results} we implement the qualitative techniques,
namely Poincar\'e maps, frequency analysis and Fast Lyapunov Indicators.
The comparison between the analytical and numerical approaches, as well as some conclusions,
are given in Section~\ref{sec:conclusion}.

\section{The three--body problem with an oblate primary and SRP}\label{sec:3BP}
In this section we start by introducing the equations of motion
describing the SCR3BP with a radiating larger primary and an
oblate smaller primary (see Section~\ref{sec:model}). Then, we proceed to compute the location of
the collinear equilibrium points as a function of the parameters
$\mu$, $\beta$, $A$ (see Section~\ref{sec:equilibrium}).

\subsection{The equations of motion}\label{sec:model}
Let us denote by $\mu$ and $1-\mu$ the masses of the primaries with $\mu \in (0, 1/2]$ and let us assume the units of measure,
so that the gravitational constant is unity and the period of the primaries is equal to $2\pi$.
We consider a synodic reference system $(O,X,Y,Z)$, rotating with the angular velocity of the primaries, with the origin located at the
barycenter of the primaries and the abscissa along the primaries' axis. The position of the smaller primary is at
$(-1+\mu,0,0)$, while the larger primary is located at $(\mu,0,0)$.
Let $(X,Y,Z)$ be the coordinates of the third body in this reference frame.
Let $(P_X,P_Y,P_Z)$ be the conjugated kinetic momenta defined as
$P_{X}= \dot{X}-nY$, $P_{Y}= \dot{Y}+nX$, $P_{Z}= \dot{Z}$, where $n$ denotes the mean motion, $n=2\pi/T_{rev}$ and
$T_{rev}$ is the period of revolution.

The equations of motion are given by
\beqa{eqmotion}
\ddot{X}-2n\dot{Y}&=&\frac{\partial \Omega}{\partial X}\nonumber\\
\ddot{Y}+2n\dot{X}&=&\frac{\partial \Omega}{\partial Y}\nonumber\\
\ddot{Z}&=&\frac{\partial \Omega}{\partial Z}\ ,
\eeqa
where
$$
\Omega=\Omega(X,Y,Z)\equiv\frac{n^{2}}{2}(X^{2}+Y^{2})+\frac{q(1-\mu)}{r_{1}}+\frac{\mu}{r_{2}}\left[ 1+\frac{A}{2r_{2}^{2}}\left( 1-\frac{3Z^{2}}{r_{2}^{2}}\right) \right]
$$
with the mean motion given by $n=\sqrt{1+\frac{3}{2}A}$ (compare with Appendix A), while the distances $r_1$, $r_2$ from the primaries are given by
$r_{1}=\sqrt{(X-\mu)^{2}+Y^{2}+Z^{2}}$ and $\displaystyle r_{2}=\sqrt{(X-\mu+1)^{2}+Y^{2}+Z^{2}}$ (see, e.g., \cite{douskos,Ob}).
We have used $q=1-\beta$, where the \sl performance \rm $\beta$ of the sail is defined by
$$
\beta\equiv\frac{L_{\odot} Q}{4 \pi c \ G M_{\odot} B}\ ,
$$
where $L_{\odot} = 3.839 \times 10^{26} \ Watt$ is the Sun luminosity, $Q$ is one plus the reflectivity,
$c$ is the speed of light, $G$ is the gravitational constant,
$M_{\odot}$ the Solar mass and $B$ the mass/area ratio of the spacecraft.

Notice that equations \equ{eqmotion} are associated to the following Hamiltonian function:
\beqa{ham}
H( X, Y, Z, P_{X}, P_{Y}, P_{Z})&=&\frac{1}{2}(P_{X}^{2}+P_{Y}^{2}+P_{Z}^{2})+nYP_{X}-nXP_{Y}-\frac{q(1-\mu)}{r_{1}}\nonumber\\
 &-&\frac{\mu}{r_{2}}\left[ 1+\frac{A}{2r_{2}^{2}}\left( 1-\frac{3Z^{2}}{r_{2}^{2}}\right) \right]\ .
\eeqa

\subsection{The equilibrium points}\label{sec:equilibrium}
It is well known (see, e.g., \cite{Szebehely}) that the circular, restricted three--body problem admits the Lagrangian equilibrium points, precisely
the triangular solutions, denoted as $L_4$ and $L_5$, and the collinear solutions, denoted as $L_1$, $L_2$, $L_3$,
the latter being located on the axis joining the primaries. Due to the oblateness and to the radiation pressure effects,
out--of--plane equilibria can be found in the spatial case (see \cite{douskos}), but we will not consider such solutions in the present work.

To locate the equilibrium positions we impose that the derivatives of $\Omega$ in these points are zero, say
${{\partial\Omega}\over {\partial X}}={{\partial\Omega}\over {\partial Y}}={{\partial\Omega}\over {\partial Z}}=0$. In particular, we obtain:
\beqano
{{\partial\Omega}\over {\partial Y}}&=&n^{2}Y-\frac{q Y(1-\mu)}{(Y^{2}+Z^{2}+(X-\mu)^{2})^{\frac{3}{2}}}+\frac{15AYZ^{2}\mu}{2(Y^{2}+Z^{2}+(1+X-\mu)^{2})^{\frac{7}{2}}}\nonumber\\
&-&\frac{3AY\mu}{2(Y^{2}+Z^{2}+(1+X-\mu)^{2})^{\frac{5}{2}}}-\frac{Y\mu}{(Y^{2}+Z^{2}+(1+X-\mu)^{2})^{\frac{3}{2}}}\ ,
\eeqano
so that we have ${{\partial\Omega}\over {\partial Y}}=0$ whenever $Y=0$. In a similar way we obtain:
\beqano
{{\partial\Omega}\over {\partial Z}}&=&-\frac{qZ(1-\mu)}{(Y^{2}+Z^{2}+(X-\mu)^{2})^{\frac{3}{2}}}+\frac{15AZ^{3}\mu}{2(Y^{2}+Z^{2}+(1+X-\mu)^{2})^{\frac{7}{2}}}\nonumber\\
&-&\frac{9AZ\mu}{2(Y^{2}+Z^{2}+(1+X-\mu)^{2})^{\frac{5}{2}}}-\frac{Z\mu}{(Y^{2}+Z^{2}+(1+X-\mu)^{2})^{\frac{3}{2}}}
\eeqano
and we have ${{\partial\Omega}\over {\partial Z}}=0$ whenever $Z=0$.
As for the derivative with respect to the first component, we have:
\beqa{eq1}
{{\partial\Omega}\over {\partial X}}&=&n^{2}X-\frac{q(1-\mu)(X-\mu)}{(Y^{2}+Z^{2}+(X-\mu)^{2})^{\frac{3}{2}}}+
\frac{15AZ^{2}\mu(1+X-\mu)}{2(Y^{2}+Z^{2}+(1+X-\mu)^{2})^{\frac{7}{2}}}\nonumber\\
&-&\frac{3A(1+X-\mu)\mu}{2(Y^{2}+Z^{2}+(1+X-\mu)^{2})^{\frac{5}{2}}}-\frac{(1+X-\mu)\mu}{(Y^{2}+Z^{2}+(1+X-\mu)^{2})^{\frac{3}{2}}}\ .
\eeqa
Inserting $Y=Z=0$ in \equ{eq1} we get the equation
\beq{eqX}
n^{2}X-\frac{q(1-\mu)(X-\mu)}{|X-\mu|^{3}}-\frac{3A\mu(1+X-\mu)}{2|1+X-\mu|^{5}}-\frac{\mu(1+X-\mu)}{|1+X-\mu|^{3}}=0\ .
\eeq
Let us denote by $\gamma_{j}$ the distance between $L_{j}$ and the closer primary. Since the
collinear points $L_1$, $L_2$, $L_3$ lie in the intervals $(-1+\mu, \mu)$, $(-\infty, -1+\mu)$, $(\mu, +\infty)$,
setting $X=\gamma_{1}+\mu-1$ for $L_{1}$, $X=-\gamma_{2}+\mu-1$ for $L_{2}$, $X=\gamma_{3}+\mu$ for $L_{3}$,
the quantity $\gamma_{j}$ is found as the unique positive solution of the following generalized Euler's equations:
\beqano
&\pm& 2n^{2}\gamma_{j}^{7}+(2n^{2}\mu-6n^{2})\gamma_{j}^{6}\pm(6n^{2}-4n^{2}\mu)\gamma_{j}^{5}+(2n^{2}\mu-2q\mu-2n^{2}+2q \mp2\mu)\gamma_{j}^{4}\nonumber\\
&+&4\mu\gamma_{j}^{3}\mp(2\mu+3A\mu)\gamma_{j}^{2}+6A\mu\gamma_{j}\mp3A\mu=0
\eeqano
for $j=1,2$, where the upper sign holds for $L_1$, while the lower sign holds for $L_2$; as for $L_3$, we have that $\gamma_3$
is the solution of the following Euler's equation
\beqano
&&2n^{2}\gamma_{3}^{7}+(8n^{2}+2n^{2}\mu)\gamma_{3}^{6}+(12n^{2}+8n^{2}\mu)\gamma_{3}^{5}+(8n^{2}+2q\mu+12n^{2}\mu-2q-2\mu)\gamma_{3}^{4}\nonumber\\
&+&(2n^{2}-8q-4\mu+8n^{2}\mu+8q\mu)\gamma_{3}^{3}+(2n^{2}\mu-3A\mu+12q\mu-12q-2\mu)\gamma_{3}^{2}\nonumber\\ &+&(8q\mu-8q)\gamma_{3}+2q\mu-2q=0\ .
\eeqano
Making use of equation \equ{eqX}, we show below the dependence of the location of the collinear points as the parameters
$(\mu,A,\beta)$ are varied, where we assume that the parameters belong to the following intervals\footnote{The upper bound on $A$, say $A\leq 10^{-4}$, is definitely
large for solar system bodies, but it might apply to extrasolar planetary systems.}$^{,}$\footnote{A realistic upper bound on
the sail performance should be $\beta\leq 0.1$; however, we consider $\beta$ up to 0.5 as it is often
done in the literature, see e.g., \cite{Ari1}.}:
\beq{parint}
0<\mu\leq 0.5\ , \qquad 0\leq A\leq 10^{-4}\ , \qquad 0\leq \beta\leq 0.5\ .
\eeq

Here and in the following sections we consider three paradigmatic cases, which correspond to a spacecraft or a solar sail orbiting in the
Earth--Moon system, in the Sun--barycenter system and in the Sun--Vesta system; these three cases encompass missions to satellites, planets or asteroids and are characterized by different values of the parameters as well as by different distances of the collinear points,
as reported in Table~\ref{tab:parameters}. \\

\section{Linear stability of the collinear points and reduction of the quadratic part}\label{sec:linearstab}
To study the stability of the collinear equilibrium points and their dependency on the parameters,
we compute the linearization of the equations of motion, which requires to expand the Hamiltonian in series up to the second order.
Indeed, we will expand it directly up to the fourth order as this will be used for the reduction to the center manifold in Section~\ref{sec:center}.
Although the computation of the linear stability is rather elementary, it is mandatory before performing the center manifold reduction as in Section~\ref{sec:center}.

\subsection{Expansion of the Hamiltonian to fourth order}\label{sec:expansion}
We shift and scale the equilibrium points by making the following transformation (\cite{JM}):
\beq{shift}
X=\mp\gamma_{j}x+\mu+\alpha\ , \qquad Y=\mp\gamma_{j}y\ , \qquad Z=\gamma_{j}z\ ,
\eeq
where $\gamma_j$ denotes again the distance of $L_j$ from the closer primary;
the upper sign corresponds to $L_{1,2}$ and the lower sign to $L_3$, while
$\alpha\equiv-1+\gamma_{1}$ for $L_{1}$, $\alpha\equiv-1-\gamma_{2}$ for $L_{2}$ and $\alpha\equiv\gamma_{3}$ for $L_{3}$.
These conventions will hold hereafter.
It must be remarked that, in all cases, the transformation \equ{shift} is symplectic of parameter $\gamma_j^2$,
so that the resulting Hamiltonian must be divided by such factor. Inserting \equ{shift} in \equ{eq1}, we have:
$$
{{\partial\Omega}\over {\partial X}}=n^{2}X+\frac{\partial}{\partial X}\left[\frac{q(1-\mu)}{r_{1}}+\frac{\mu}{r_{2}}+\frac{A\mu}{2r_{2}^{3}} -\frac{3\gamma_1^2 z^2 A\mu}{2r_{2}^{5}}\right]\ ;
$$
therefore, we obtain the equation of motion:
$$
\mp\ddot{x}\pm2n\dot{y}\pm n^{2}x= \frac{n^{2}(\mu+\alpha)}{\gamma_{j}}\mp\frac{1}
{\gamma_{j}^{2}}\frac{\partial}{\partial x}\left[\frac{q(1-\mu)}{r_{1}}+
\frac{\mu}{r_{2}}+\frac{A\mu}{2r_{2}^{3}} -\frac{3\gamma_j^2 z^2A\mu}{2r_{2}^{5}}\right]\ ,
$$
where $j=1,2,3$. We remark that the inverse of the distances are transformed as
$$
\frac{1}{r_{1}}=\frac{1}{\gamma_{j}\sqrt{(x\mp\frac{\alpha}{\gamma_{j}})^{2}+y^{2}+z^{2}}}\ ,\qquad \frac{1}{r_{2}}=\frac{1}{\gamma_{j}\sqrt{(-x\pm\frac{\alpha+1}{\gamma_{j}})^{2}+y^{2}+z^{2}}}\ .
$$
In a similar way, for the other two components we obtain:
\beqano
\mp\gamma_{j}\ddot{y}\mp2n\gamma_{j}\dot{x}&=&\mp n^{2}(\gamma_{j} y)\mp\frac{1}{\gamma_{j}}\frac{\partial}{\partial y}
\left[\frac{q(1-\mu)}{r_{1}}+\frac{\mu}{r_{2}}+\frac{A\mu}{2r_{2}^{3}} -\frac{3\gamma_1^2 z^2 A\mu}{2r_{2}^{5}}\right]\nonumber\\
\gamma_{j}\ddot{z}&=&\frac{1}{\gamma_j}\frac{\partial}{\partial z}\left[\frac{q(1-\mu)}{r_{1}}+\frac{\mu}{r_{2}}+
\frac{A\mu}{2r_{2}^{3}} -\frac{3 \gamma_1^2 z^2 A\mu}{2r_{2}^{5}}\right]\ .
\eeqano
Expanding in Taylor series $1/r_1$, $1/r_2$ and their powers, we compute the linearized equations as
\beqano
\ddot{x}-2n\dot{y}-(n^{2}-2a)x&=& 0\nonumber\\
\ddot{y}+2n\dot{x}+(-n^{2}+b)y&=& 0\nonumber\\
\ddot{z}+cz&=& 0\ ,
\eeqano
where the quantities $a$, $b$, $c$ take different values according to the considered equilibrium point.
Since in the following sections we shall analyze only $L_1$ and $L_2$, we provide the explicit values for such positions.
Precisely, one has:
\beqano
a&=&\frac{q(1-\mu)}{\alpha^{3}}\mp \frac{\mu}{(1+\alpha)^{3}}\mp \frac{3A\mu}{(1+\alpha)^{5}}\nonumber\\
b&=&-\frac{q(1-\mu)}{\alpha^{3}}\pm \frac{\mu}{(1+\alpha)^{3}}\pm\frac{3A\mu}{2(1+\alpha)^{5}}\nonumber\\
c&=&-\frac{q(1-\mu)}{\alpha^{3}}\pm \frac{\mu}{(1+\alpha)^{3}}\pm \frac{9A\mu}{2(1+\alpha)^{5}}\ ,
\eeqano
where the upper sign holds for $L_1$ and the lower sign for $L_2$.

Setting for shortness $\Delta \equiv \frac{3 A \mu}{2|1+\alpha|^5}$, we can write $a$ and $c$ in terms of $b$ as
\beq{abcsmart}
a=-(b+\Delta)\ ,\qquad c=b+2\Delta\ .
\eeq
It results that $\Delta>0$ and $b>0$ for all $A>0$, $\mu>0$, $\beta>0$.\\

Finally, the complete equations of motion in the new variables can be written as
\beqa{em}
\ddot{x}-2n\dot{y}-(n^{2}-2a)x&=& \frac{1}{\gamma_{j}^{2}}\frac{\partial}{\partial x}\sum_{n\geq 3}H_{n}\nonumber\\
\ddot{y}+2n\dot{x}+(-n^{2}+b)y&=& \frac{1}{\gamma_{j}^{2}}\frac{\partial}{\partial y}\sum_{n\geq 3}H_{n}\nonumber\\
\ddot{z}+cz&=& \frac{1}{\gamma_{j}^{2}}\frac{\partial}{\partial z}\sum_{n\geq 3}H_{n}\ ,
\eeqa
where $H_{n}$ are suitable polynomials of degree $n$.
Defining the momenta as $p_{x}=\dot{x}-ny$, $p_{y}=\dot{y}+nx$, $p_{z}=\dot{z}$ and using \equ{em},
the Hamiltonian function \equ{ham} becomes
\beq{hamexp}
H(x, y, z, p_{x}, p_{y}, p_{z})=\frac{1}{2}(p_{x}^{2}+p_{y}^{2}+p_{z}^{2})+ nyp_{x}-nxp_{y}+ax^{2}+\frac{1}{2}by^{2}+\frac{1}{2}cz^{2}-\frac{1}{\gamma_{j}^{2}} \sum_{n\geq 3}H_{n}\ .
\eeq

\subsection{Reduction of the quadratic terms and linear stability}\label{sec:reduction}
The quadratic part of the Hamiltonian \equ{hamexp} is given by
$$
H_{2}(x, y, z, p_{x}, p_{y}, p_{z})=\frac{1}{2}(p_{x}^{2}+p_{y}^{2})+ nyp_{x}-nxp_{y}+\frac{1}{2}p_{z}^{2}+ ax^{2}+\frac{1}{2}by^{2}+\frac{1}{2}cz^{2}\ .
$$
By (\ref{abcsmart}) it results that $c>0$ for each of the equilibrium points.
This implies that the vertical direction is described by a harmonic oscillator with frequency $\omega_{2}=\sqrt{c}$, namely:
$$
 \dot{p_{z}}=-\frac{\partial H}{\partial z}=-cz ,\ \ \ \ \dot{z}=\frac{\partial H}{\partial p_{z}}=p_{z}\ .
$$
As for the planar directions, following \cite{JM} the quadratic part of the Hamiltonian, that we keep denoting
as $H_2$, is given by
$$
 H_{2}(x, y, p_{x}, p_{y})=\frac{1}{2}(p_{x}^{2}+p_{y}^{2})+ nyp_{x}-nxp_{y}+ ax^{2}+\frac{1}{2}by^{2}\ .
$$
Denoting by $J$ the symplectic matrix, the equations of motion are given by
$$
\left(%
\begin{array}{c}
\dot{x}\\
\dot{y}\\
\dot{p_{x}}\\
\dot{p_{y}}\\
\end{array}
\right)
=J Hess(H_{2})
\left(\begin{array}{c}
x\\
y\\
p_{x}\\
p_{y}\\
\end{array}\right)
=\left(%
\begin{array}{cccc}
0 & n & 1 & 0\\
-n & 0 & 0 & 1\\
-2a & 0 & 0 & n\\
0 & -b & -n & 0\\
 \end{array}%
\right)
\left(%
\begin{array}{c}
x\\
y\\
p_{x}\\
p_{y}\\
\end{array}
\right)\ .
$$
Let us define the matrix $M$ as
\beq{Mdef}
M\equiv J\ Hess(H_{2})\ ;
\eeq
the associated characteristic polynomial is given by
\beq{plam}
p(\lambda)= \lambda^{4}+(2n^{2}+2a+b)\lambda^{2}+(n^{4}-2an^{2}-bn^{2}+2ab)\ .
\eeq
Setting $\eta=\lambda^{2}$, the roots of the polynomial \equ{plam} are given by
\beqano
\eta_{1}&=&\frac{-2n^{2}-2a-b-\sqrt{16an^{2}+4a^{2}+8bn^{2}-4ab+b^{2}}}{2}\nonumber\\
\eta_{2}&=&\frac{-2n^{2}-2a-b+\sqrt{16an^{2}+4a^{2}+8bn^{2}-4ab+b^{2}}}{2}\ .
\eeqano
In order to study the stability of the collinear equilibrium points we have to establish the domains
in which $\lambda_{1,3}=\pm\sqrt{\eta_1}$ and $\lambda_{2,4}=\pm\sqrt{\eta_2}$ are real, complex or imaginary.
In particular, a given equilibrium point $L_j$ will be linearly stable, if $\eta_1$ and $\eta_2$ are purely imaginary,
while it is unstable elsewhere. We are interested to the case in which the collinear points are of the type
\emph{saddle}$\times$\emph{center}$\times$\emph{center}, which occurs whenever
$\eta_1<0$ and $\eta_2>0$. This is equivalent to require that the following inequalities are satisfied:
\begin{equation}
\label{con}
\left\{\begin{array}{ll}
16an^{2}+4a^{2}+8bn^{2}-4ab+b^{2}\geq0\nonumber\\
-2n^{2}-2a-b+\sqrt{16an^{2}+4a^{2}+8bn^{2}-4ab+b^{2}}>0\nonumber\\
-2n^{2}-2a-b-\sqrt{16an^{2}+4a^{2}+8bn^{2}-4ab+b^{2}}<0\ .\nonumber\\
\end{array}\right.
\end{equation}
Using (\ref{abcsmart}), the inequalities (\ref{con}) become:
\beq{eqbb}
\left\{
\begin{array}{ll}
9 b^2 + 4 \Delta (-4 n^2 + \Delta) + b (-8 n^2 + 12 \Delta)\geq0\\
b - 2 n^2 + 2 \Delta+\sqrt{9 b^2 + 4 \Delta (-4 n^2 + \Delta) + b (-8 n^2 + 12 \Delta)}>0\\
b - 2 n^2 + 2 \Delta-\sqrt{9 b^2 + 4 \Delta (-4 n^2 + \Delta) + b (-8 n^2 + 12 \Delta)}<0\ .\\
\end{array}\right.
\eeq
The first condition in \equ{eqbb} is satisfied whenever
\beq{eqb}
b>\frac{2}{9} (2 n^2 - 3 \Delta + 2 \sqrt{n^4 + 6 n^2 \Delta})\ .
\eeq
Given that the function at the right hand side of \equ{eqb} reaches its maximum $n^2$ for $\Delta=\frac{n^2}{2}$,
the inequality \equ{eqb} is verified if $b>n^2=\sqrt{1+\frac{3}{2}A}$.
The second and third conditions in \equ{eqbb} can be reduced to verify just that
$$
8 b^2 +b(8\Delta- 4 n^2) - 8 n^2 \Delta- 4 n^4 >0\ ,
$$
which is satisfied again if $b>n^2$, which holds for all values of the oblateness, the solar radiation pressure and
the mass parameter in the intervals defined in \equ{parint}. This concludes the discussion of the stability character of the
collinear points, including the solar radiation pressure and the oblateness of one of the primaries.

\section{Center manifold reduction}\label{sec:center}
Due to the $saddle\times center \times center$ character of the collinear equilibrium points (see Section~\ref{sec:linearstab}), we proceed to
perform the reduction to the center manifold. The adopted procedure is a straightforward extension of that used in \cite{JM},
provided that we include the necessary modifications to consider the effects of
the oblateness and the solar radiation pressure. For completeness, we report here the main steps to treat the Hamiltonian
\equ{ham} (see Appendix C for more details).
We stress that after removing the hyperbolic direction, we will be able to perform a qualitative analysis based on
Poincar\'e sections, frequency analysis and Fast Lyapunov Indicators as it is done in Section~\ref{sec:results}.\\

Taking into account that $\eta_{1}<0$ and $\eta_{2}>0$,
let $\omega_{1}=\sqrt{-\eta_{1}}$ and $\lambda_{1}=\sqrt{\eta_{2}}$; we look for a change of variables, so that
we reach a simpler form of the Hamiltonian. As described in detail in Appendix B,
this is obtained by computing the eigenvalues and eigenvectors of $M$ in \equ{Mdef},
which provide a transformation allowing to get a Hamiltonian with the following quadratic part
(with a slight abuse we keep the same notation for all variables):
\beq{eq5}
H_{2}(x, y, z, p_{x}, p_{y}, p_{z})=\lambda_{1}xp_{x}+\frac{\omega_{1}}{2}(p_{y}^{2}+y^{2})+\frac{\omega_{2}}{2}(p_{z}^{2}+z^{2})\ ,
\eeq
where
\beqano
\lambda_{1}&=&\sqrt{\frac{-2n^{2}-2a-b+\sqrt{16an^{2}+4a^{2}+8bn^{2}-4ab+b^{2}}}{2}}\nonumber\\
\omega_{1}&=&\sqrt{-\frac{-2n^{2}-2a-b-\sqrt{16an^{2}+4a^{2}+8bn^{2}-4ab+b^{2}}}{2}}\nonumber\\
\omega_{2}&=&\sqrt{c}\ .
\eeqano
In analogy to \cite{JM}, we introduce the following complex transformation
\beqano
x&=&q_{1}, \ \ \ \ \ y=\frac{q_{2}+ip_{2}}{\sqrt{2}}, \ \ \ \ z=\frac{q_{3}+ip_{3}}{\sqrt{2}}\nonumber\\
p_{x}&=&p_{1}, \ \ \ \  p_{y}=\frac{iq_{2}+p_{2}}{\sqrt{2}},\ \ \
p_{z}=\frac{iq_{3}+p_{3}}{\sqrt{2}}\ , \eeqano which provides a
complex expression for the Hamiltonian; we report here the
quadratic part of the Hamiltonian, which takes the form: \beq{eq6}
H_{2}(q_{1}, q_{2}, q_{3}, p_{1}, p_{2},
p_{3})=\lambda_{1}q_{1}p_{1}+i\omega_{1}q_{2}p_{2}+i\omega_{2}q_{3}p_{3}\ .
\eeq
Beside the quadratic part \equ{eq6}, we need to compute the nonlinear terms.
Straightforward but tedious computations, performed by means of the $\texttt{Mathematica}^\copyright$ algebraic manipulator,
allow us to find the expressions of $H_{3}$ and $H_{4}$ in complex variables. Afterwards we shall implement
a Lie series transformation to obtain the reduction to the center manifold as described in the following section.

\subsection{Reduction to the center manifold}\label{sec:CM}
The reduction to the center manifold is obtained by making suitable changes of variables
by using Lie series (see Appendix C). Indeed, we conjugate the SCR3BP to a Hamiltonian of the form
$$
 H(q,p)= \lambda_1 q_{1}p_{1}+\textit{i}\omega_{2}q_{2}p_{2}+\textit{i}\omega_{3}q_{3}p_{3}+ \sum_{n\geq 3}H_{n}(q,p)\ ,
$$
for suitable coordinates $(q_1,q_2,q_3)$ and momenta $(p_1,p_2,p_3)$, where the quadratic part has been obtained in
\equ{eq6} and $H_{n}$ denotes a homogeneous polynomial of degree $n$.
In the linear approximation the center manifold is obtained by imposing $q_{1}=p_{1}=0$, since the
hyperbolicity pertains to such variables. Then, we will require that $\dot{q_{1}}(0)=\dot{p_{1}}(0)=0$ for $q_{1}(0)=p_{1}(0)=0$,
so that we obtain $q_{1}(t)=p_{1}(t)=0$ for any time, due to the autonomous character of the problem. Taking into account Hamilton's equations
$$
\dot{q_{i}}=H_{p_{i}}\ , \qquad \dot{p_{i}}=-H_{q_{i}}\ ,
$$
this requirement is satisfied whenever in the series expansion of the
Hamiltonian all monomials of the type $h_{ij}q^{i}p^{j}$ with $i_{1}\not= j_{1}$ are such that $h_{ij}=0$, being $i=(i_{1},i_{2},i_{3})$ and  $j=(j_{1},j_{2},j_{3})$. In this way we obtain a Hamiltonian of the form $H(q,p)= H_{N}(q,p)+R_{N}(q,p)$, where
$H_{N}(q,p)$ is a polynomial of degree $N$ in $(q,p)$ without terms depending on the product $q_1p_1$, while
$R_{N}(q,p)$ is a reminder of order $N+1$. We refer to Appendix C for the description of
a procedure based on Lie series to determine explicitly the
required canonical transformation. Let us denote by $(y,z,p_y,p_z)$ the normalized variables; the final expression of the Hamiltonian
reduced to the center manifold has the following form:
\beq{hamCM}
\tilde H(y,z,p_y,p_z)={\omega_1\over 2}(p_y^2+y^2)+{\omega_2\over 2}(p_z^2+z^2)+\tilde H_3(y,z,p_y,p_z)+\tilde H_4(y,z,p_y,p_z)\ ,
\eeq
where
$\tilde H_3$, $\tilde H_4$ denote homogeneous polynomials of
degree, respectively, 3 and 4. The frequencies $\omega_1$, $\omega_2$, as well as those of the higher order terms,
depend on the choice of the parameters and will be specified in each concrete case.

\section{Analytical estimates of the bifurcation thresholds}\label{sec:analytical}
Analytical results providing an estimate for the value of the thresholds at which the bifurcation of halo
orbits takes place have been presented in \cite{CPS}. The result is briefly summarized as follows. After the reduction to the
center manifold, a normal form is computed around the synchronous resonance. The resulting normal form admits a
first integral, related to the action variables of the harmonic oscillator (i.e., the quadratic part in \equ{hamCM}).
A detuning measuring the displacement around the synchronous resonance is introduced; assuming that the detuning is small, one can
compute the bifurcation threshold at different orders in the powers series expansion in the detuning.
In \cite{CPS} the computation at first and second order has been performed. Here, we extend the method of
\cite{CPS} by computing the thresholds for the model including oblateness and solar radiation pressure.
As we will see in Section~\ref{sec:results}, the case in which the parameter $\beta$ is different from zero allows one
to find several bifurcations at relatively low energy levels. We anticipate that the analytical estimates computed in this section will
agree with the numerical values of Section~\ref{sec:results} (see Table~\ref{tab:bif}).\\

Let us write \equ{hamCM} in complex form, implementing the following change of coordinates:
\beqano
Q_2&=&{{\sqrt{2}}\over {2i}} (p_2+iq_2)\ ,\qquad Q_3={{\sqrt{2}}\over {2i}} (p_3+iq_3)\ ,\nonumber\\
P_2&=&{{\sqrt{2}}\over {2}} (p_2-iq_2)\ ,\qquad P_3={{\sqrt{2}}\over {2}} (p_3-iq_3)\ .
\eeqano
Next, we introduce action--angle variables $(I_2,I_3,\theta_2,\theta_3)$ for the quadratic part of the Hamiltonian by means
of the change of coordinates:
\beqano
Q_2&=&-i\sqrt{I_2}\ e^{i\theta_2}\ ,\qquad Q_3=-i\sqrt{I_3}\ e^{i\theta_3}\ ,\nonumber\\
P_2&=&\sqrt{I_2}\ e^{-i\theta_2}\ ,\qquad\ \  P_3=\sqrt{I_3}\ e^{-i\theta_3}\ ,
\eeqano
so that we obtain a Hamiltonian of the form
$$
H(\theta_2,\theta_3,I_2,I_3)=\omega_2 I_2+\omega_3 I_3+\tilde H_3(\theta_2,\theta_3,I_2,I_3)+\tilde H_4(\theta_2,\theta_3,I_2,I_3)\ .
$$
To investigate the appearance of the resonant periodic orbits we compute a normal form in the neighborhood of the synchronous resonance $\omega_2=\omega_3$.
The explicit computation shows that the first non--trivial order is given by $\tilde H_4$, since the third degree term
$\tilde H_3$ does not contain resonant terms. We are thus led to a resonant normal form given by
\beq{NF}
H_{NF}(\theta_2,\theta_3,I_2,I_3)=\omega_2 I_2+\omega_3 I_3+\Big[a_{20}I_2^2+a_{02}I_3^2+I_2I_3(a_{11}+
2b_{11}\cos(2\theta_2-2\theta_3)\Big]\ ,
\eeq
where the coefficients $a_{20}$, $a_{02}$, $a_{11}$, $b_{11}$ are evaluated explicitly in Table~\ref{tab:ab}.

The dynamics is determined by the \sl normal modes \rm $I_k = {\rm const.}$, $k=2,3$ and by the periodic orbits in general position,  related to the resonance. The normal modes always exist and at low energies are both stable: $I_2 = {\rm const.}, \; I_3=0$ gives rise to the \sl planar Lyapunov orbit; \rm $I_2 = 0, \; I_3=  {\rm const.}$ gives rise to the \sl vertical Lyapunov orbit. \rm When stable, these orbits are surrounded by families of \sl Lissajous \rm tori. \rm The resonant families may appear as bifurcation from the normal mode at some given energy threshold and are determined by the condition that the frequency
of the normal mode is equal to that of its normal perturbation. The normal modes can be again stable through a second bifurcation; a concrete example of a
second bifurcation for the Earth--Moon case is given in \cite{GM} (we refer to \cite{CPS,Pucacco,Pucacco14,Pucacco14bis} for
further details). To investigate the possible sequences of bifurcations we proceed as follows.

\vskip.1in

\begin{table}
\caption{Coefficients of the normal form \equ{NF} around $L_1$ for the Earth--Moon, Sun--barycenter, Sun--Vesta systems
with and without solar radiation pressure.}\label{tab:ab}
\center
\tabcolsep=3mm
\renewcommand\arraystretch{1.2}
\begin{tabular}{|c|c|c|c|c|}
\hline
& $a_{20}$& $a_{02}$ & $a_{11}$ & $b_{11}$\\
\hline
\hline
Earth--Moon & $0.162109 $ & $0.144891 $ & $0.0726274 $& $0.116537 $\\
\hline
Sun--barycenter $\beta=0$ & $0.0989667 $ & $0.08098 $ & $-0.0256235 $& $0.102138 $\\
\hline
Sun--barycenter $\beta=10^{-2}$ & $0.136123 $ & $0.108011 $ & $0.0189407 $& $0.11106 $\\
\hline
Sun--Vesta $\beta=0$ & $0.0957347$ & $0.0777063$ & $-0.0304854$& $0.101319$\\
\hline
Sun--Vesta $\beta=10^{-2}$ & $0.0157472 $ & $0.00203253 $ & $4.11966\ 10^{-7} $& $0.00533371 $\\
\hline
\end{tabular}
\end{table}

\vskip.1in

From Hamilton's equations associated to \equ{NF}, it is readily seen that $\dot I_2+\dot I_3=0$. This remark leads
to introduce the following change of variables (\cite{CPS}):
\beqa{ER}
\E&=&I_2+I_3\ ,\qquad \R=I_2\ ,\nonumber\\
\nu&=&\theta_3\ ,\qquad\qquad\  \psi=\theta_2-\theta_3\ .
\eeqa
Moreover, following \cite{CPS} let us introduce the \sl detuning \rm $\delta$ as
$$
\delta=\omega_2-\omega_3\ ,
$$
which measures the displacement from the synchronous resonance. Using \equ{ER} the Hamiltonian \equ{NF} becomes
$$
H_{new}(\E,\R,\nu,\psi)=\E+\tilde\delta \R+a\R^2+b\E^2+c\E\R+d(\R^2-\E\R)\ \cos 2\psi\ ,
$$
where $\tilde\delta=\delta/\omega_3$, $a=(a_{20}+a_{02}-a_{11})/\omega_3$, $b=a_{02}/\omega_3$,
$c=(a_{11}-2a_{02})/\omega_3$, $d=-2b_{11}/\omega_3$.

It has been shown in \cite{CPS} that the equilibria associated to Hamilton's equations of $H_{new}$ can be classified as
\sl inclined \rm (or `anti-halo') when $\psi=0$ or $\psi=\pi$, and \sl loop \rm (or `halo'), when $\psi=\pm \pi/2$.
In view of the reflection symmetries, each case actually corresponds to a double family. They
exist at the following level values of the integral of motion $\E$:
\beq{inclined}
\E_{iy}={{\delta\omega_2^2}\over {-a_{11}+2(a_{20}-b_{11})}}\ , \qquad
\E_{iz}={{\delta\omega_2^2}\over {-2a_{02}+a_{11}+2b_{11}}}\ ,
\eeq
for the inclined families and
\beq{loop}
\E_{\ell y}={{\delta\omega_2^2}\over {-a_{11}+2(a_{20}+b_{11})}}\ , \qquad
\E_{\ell z}={{\delta\omega_2^2}\over {-2a_{02}-a_{11}+2b_{11}}}\ ,
\eeq
for the loop families. These values correspond to a first order computation in the detuning;
refined (but more complicated) values at second order have been found analytically in \cite{CPS}.

The physical interpretation of the thresholds in \equ{inclined}, \equ{loop} is the following.
The quantity $\E_{\ell y}$ determines the bifurcation of the halo families from the planar Lyapunov orbit, when this becomes
unstable. At $\E_{iy}$ the planar Lyapunov orbit turns back to being stable and one observes the
bifurcation of the (unstable) anti-halo orbits. Finally, the two unstable families which have been formed collapse on the
vertical at $\E_{iz}$. The disappearance of the halo at $\E_{\ell z}$ never occurs in all cases we have investigated.

This sequence of bifurcations will be clearly shown in the case of the asteroid Vesta under the
effect of solar radiation pressure (see Section~\ref{sec:results}).

\section{Qualitative analysis of the bifurcation values}\label{sec:results}
On the basis of the center manifold reduction obtained in Section~\ref{sec:center}, we implement some numerical techniques
which allow us to characterize the dynamics and, in particular, to distinguish between the different types of orbits,
precisely planar Lyapunov and halo orbits, the latter ones obtained at specific levels of the energy at which the
bifurcation takes place.
As concrete models, we consider three paradigmatic cases, characterized by different values of the mass ratio: a
relatively high value as in the Earth--Moon system, an intermediate mass ratio as for the Sun--barycenter system
(between the barycenter of the Earth--Moon system and the Sun), and
a low value as in the Sun--Vesta system.

\subsection{Poincar\'e section}\label{sec:poinc}
To get a qualitative description of the dynamics in the center manifold we start by computing a Poincar\'e section as
follows. We set $z=0$ and we fix an energy level $H=h_{0}$, from which we compute the initial value of $p_{z}$ choosing the
solution with $p_{z}>0$. The Poincar\'e section is then shown in the plane $(y,p_y)$; we will see that, as the energy increases and
exceeds a specific energy value, halo orbits arise from bifurcations of planar Lyapunov periodic orbits.

\subsection{Frequency analysis}\label{sec:freq}
This technique consists in studying the behavior of the \sl frequency map \rm (\cite{Laskar,LFC}),
which is obtained computing the variation of the absolute value of the ratio of the frequencies,
say $\omega_r=|\omega_y/\omega_z|$, as a function of the initial values of the
action variables, whereas the initial conditions of the angles can be set to zero (see, e.g., \cite{CFL}, see also \cite{CF}).

The frequency analysis has the advantage to be computationally fast and it allows us to obtain a complementary investigation
of the occurrence of halo orbits. Precisely, we proceed as follows. Concerning the initial conditions,
we fix as starting values $z=0$ and $p_y=p_y^0$, we scan over initial values for $y$ in a given interval
and for an assigned energy level $H=h_0$, we compute the corresponding value of $p_z$.
We find convenient to avoid using Cartesian variables, and we rather transform to action-angle variables
for the quadratic part of \equ{hamCM}. Thus, we introduce harmonic oscillator actions $(J_y,J_z)$ defined through
the expressions
\beqano
p_y&=&\sqrt{2J_y}\cos\theta_y\ ,\qquad y=\sqrt{2J_y}\sin\theta_y\ ,\nonumber\\
p_z&=&\sqrt{2J_z}\cos\theta_z\ ,\qquad z=\sqrt{2J_z}\sin\theta_z\ ,
\eeqano
where $(\theta_y,\theta_z)$ denote the conjugated angle variables. Next we perform a first order perturbation theory
by averaging over the angle variables to obtain a normalized Hamiltonian, whose derivative provides an expression
for the frequencies associated to the given initial conditions. Finally, we back--transform to the variables $(J_y,J_z)$
to get a frequency vector $(\omega_y,\omega_z)$ associated to the previous initial data. The frequency map is
obtained by computing the variation of $\omega_r=|\omega_y/\omega_z|$ as the initial condition is varied.

\subsection{Fast Lyapunov Indicator}\label{sec:FLI}
In order to investigate the stability of the dynamics in the center manifold, we compute a quantity called the Fast Lyapunov Indicator
(hereafter FLI), which is determined as the value of the largest Lyapunov characteristic exponent
at a fixed time (see \cite{froes}). By comparing the values of the FLIs as the initial conditions or suitable parameters are varied, one obtains
an indication of the dynamical character of the orbits (precisely Lyapunov or halo) as well as of their stability.
The explicit computation of the FLI proceeds as follows. Let ${\underline \xi}=(y,z,p_y,p_z)$,
let the vector field associated to \equ{hamCM} be denoted as
$$
\dot{\underline \xi}={\underline f}({\underline \xi})\ ,\qquad {\underline{\xi}}\in{\bf R}^4\ ,
$$
and let the corresponding variational equations be
$$
\dot{\underline \eta}=\Big({{\partial \underline{f}(\underline{\xi})} \over {\partial
\underline{\xi}}}\Big)\ {\underline \eta}\ ,\qquad {\underline{\eta}}\in{\bf R}^4\ .
$$
Having fixed an initial condition $\underline{\xi}(0) \in {\bf R}^{4}$,
$\underline{\eta}(0) \in {\bf R}^{4}$, the FLI at a given time $T\geq 0$ is obtained by the expression
$$
{\rm FLI}(\underline{\xi}(0),  \underline{\eta}(0), T) \equiv \sup _{0 < t\leq T}
 \log || \underline{\eta}(t)||\ ,
$$
where $\|\cdot\|$ denotes the Euclidean norm.

\subsection{Applications}\label{sec:applications}
We proceed to implement the techniques described in Section~\ref{sec:poinc}, \ref{sec:freq}, \ref{sec:FLI} to
the concrete samples provided by the Earth--Moon, Sun--barycenter, Sun--Vesta systems. The computations have been
performed using $\texttt{Mathematica}^\copyright$ as well as developing dedicated programs in a general-purpose programming language.
The parameters associated to these three samples are listed in Table~\ref{tab:parameters}.

\vskip.1in

\begin{table}
\caption{Main parameters and location (in normalized units) of the collinear points for the Earth--Moon, Sun--barycenter, Sun--Vesta systems.}\label{tab:parameters}
\center
\tabcolsep=3mm
\renewcommand\arraystretch{1.2}
\begin{tabular}{|c|c|c|c|}
\hline
\hline
& Earth--Moon& Sun--barycenter & Sun--Vesta \\
\hline
\hline
$\mu$ & $1.2154\cdot 10^{-2}$ & $3.040423\cdot 10^{-6}$ & $1.3574\cdot 10^{-10}$\\
\hline
$J_{2}$ & $2.034\cdot 10^{-4}$ & $1081\cdot 10^{-6}$ & $0.0812232$ \\
\hline
$A$ & $4.15559\cdot 10^{-9}$ & $1.96782\cdot 10^{-12}$ & $4.54776\cdot 10^{-14}$\\
\hline
$\beta$ & 0 & $10^{-2}$ & $10^{-2}$\\
\hline
$L_1$ & -0.836898 & -0.988731 & -0.996651\\
\hline
$L_2$ & -1.1557 & -1.00908 & -1.000115\\
\hline
$L_3$ & 1.00506 & 0.996657 & 0.996655 \\
\hline
\end{tabular}
\end{table}

\vskip.1in

The values of the quantities introduced in Section~\ref{sec:center}, needed for the center manifold
reduction, are listed in Table~\ref{tab:parcm}.

\vskip.1in

\begin{table}
\caption{Quantities for the center manifold reduction associated to the Earth--Moon, Sun--barycenter, Sun--Vesta systems.}\label{tab:parcm}
\center
\tabcolsep=3mm
\renewcommand\arraystretch{1.2}
\begin{tabular}{|c|c|c|c|}
\hline
\hline
& Earth--Moon& Sun--barycenter & Sun--Vesta \\
\hline
\hline
$\gamma_{1}$ & $0.150948$ & $0.01127$ & $0.00334854$ \\
\hline
$\alpha (\gamma_1)$& $-0.836898$ & -$0.98873$ & -0.996651 \\
\hline
$a$ & $-5.14772$ & -3.15056 & -1.00363 \\
\hline
$b$ & $5.14772$ & 3.15056 & 1.00363\\
\hline
$c$ & $5.14772$ & 3.15056 & 1.00363\\
\hline
$\lambda_{1}$ &2.9321 & 2.13994 & 0.10407 \\
\hline
$\omega_{1}$ & 2.33441 & 1.85169 & 1.0036\\
\hline
$\omega_{2}$ & 2.26886 & 1.77498 & 1.00181\\
\hline
$s_{1}$ & 14.9084 & 9.6584 & 0.79682 \\
\hline
$s_{2}$ & 23.4324 & 12.6138 & 2.02523 \\
\hline
\end{tabular}
\end{table}

\vskip.1in

\begin{figure}
\centering
\includegraphics[width=0.28\linewidth]{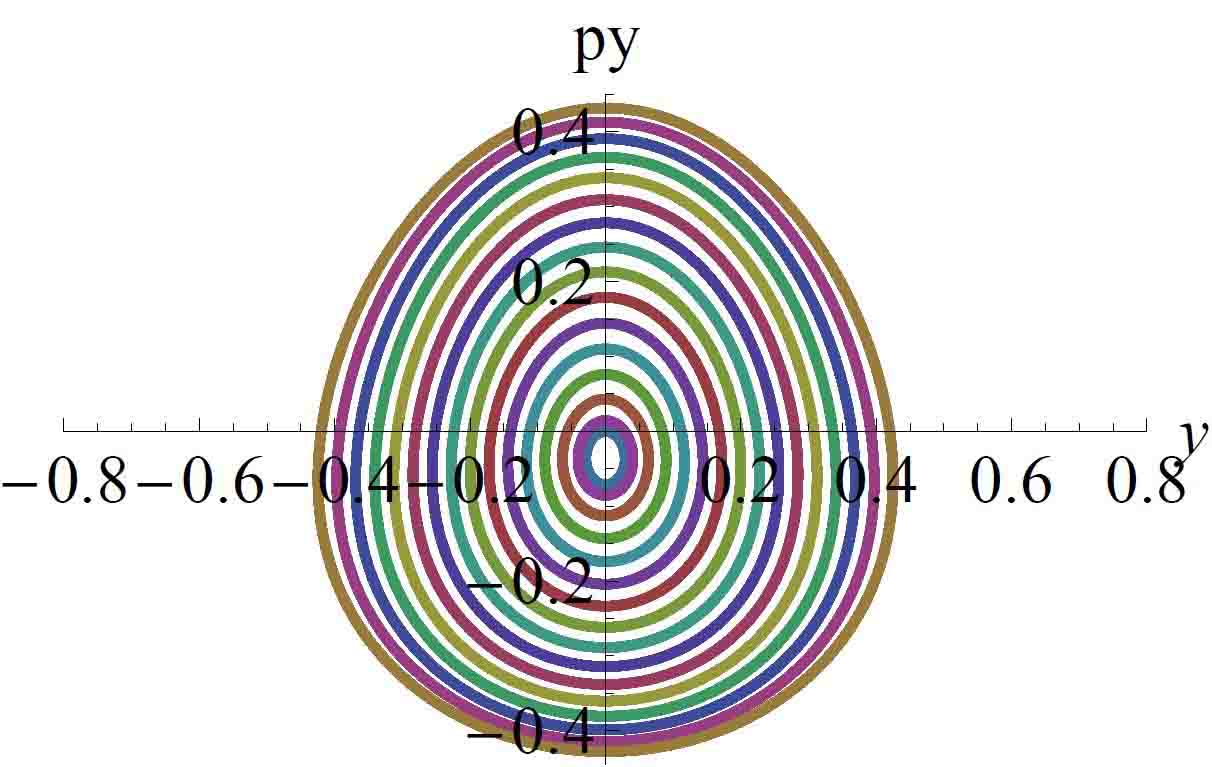}
\includegraphics[width=0.28\linewidth]{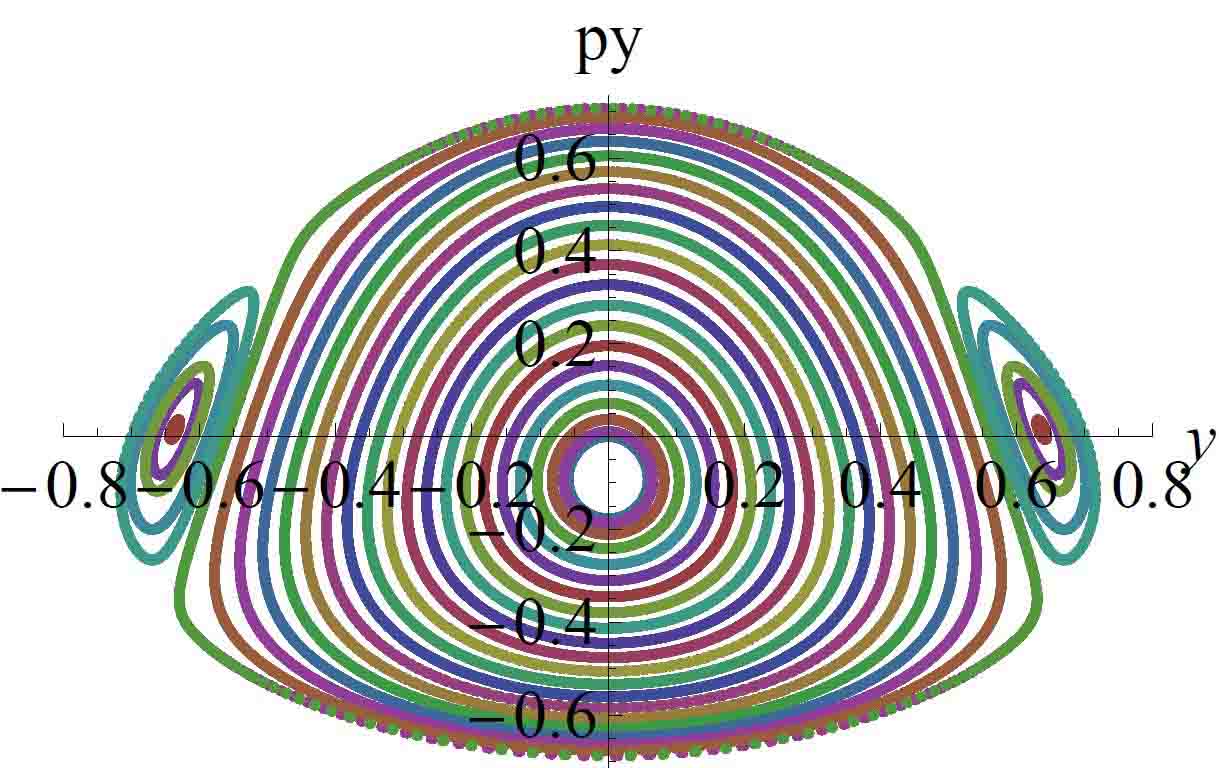}
\includegraphics[width=0.28\linewidth]{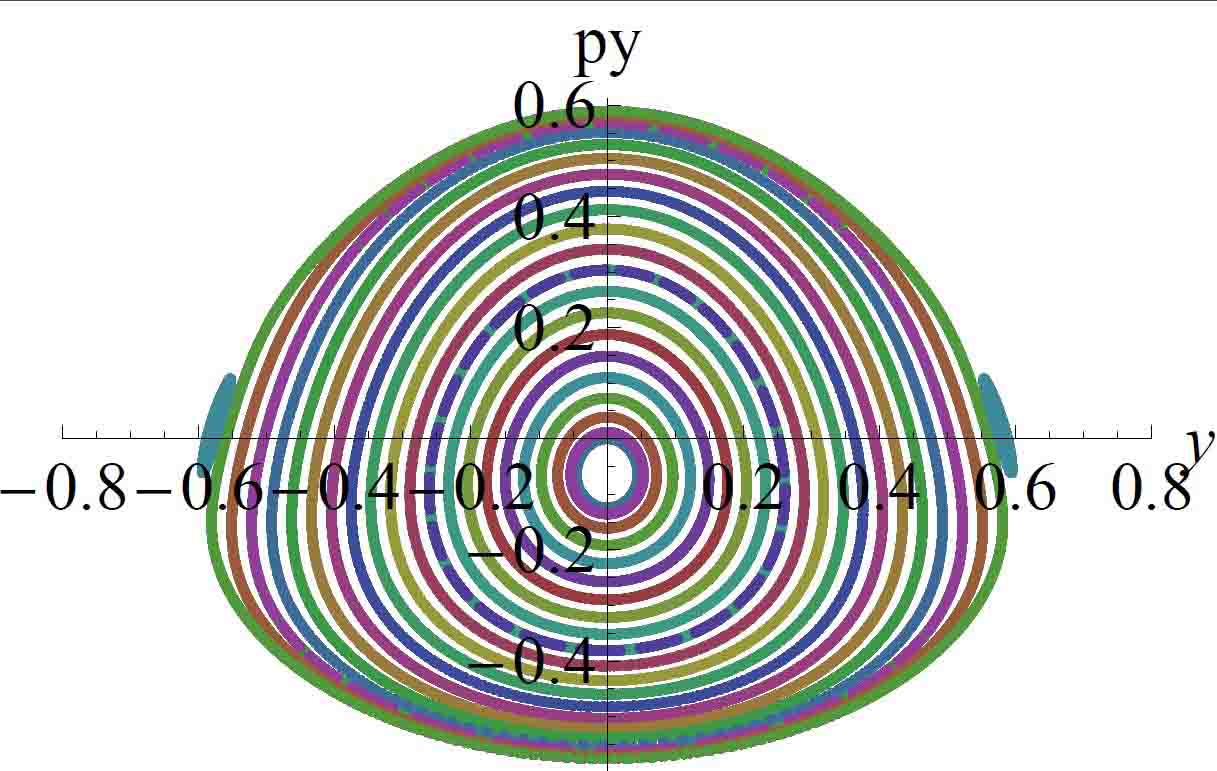}
\caption{Poincar\'e sections of the center manifold associated to $L_1$ of the Sun--Vesta system on the plane $(y,p_{y})$
without solar radiation pressure (i.e., $\beta=0$) and no oblateness (i.e., $A=0$);
different values of the energy are taken into account: $h=0.2$ (left panel), 0.35 (middle), 0.5 (right).}\label{fig:Vestasenzaq}
\end{figure}

We analyze in detail the Sun--Vesta case, which presents several interesting features as the different parameters are varied.
We start by computing the Poincar\'e surfaces of section of the center manifold associated to $L_1$.
We report in Figure~\ref{fig:Vestasenzaq} the Poincar\'e sections in the plane $(y,p_{y})$ with only the
gravitational effect, namely with $\beta=0$ (no solar radiation pressure) and $A=0$ (no oblateness).
The maps show that the appearance of halo orbits occurs for an energy value approximately equal
to $h=0.35$ (more accurate computations will provide the bifurcation vale $h=0.3341$, compare
with Table~\ref{tab:bif}).
For higher levels of the energy, the amplitude of the halo orbits increases as shown
by the map at $h=0.5$.

\begin{figure}
\centering%
\includegraphics[width=0.35\linewidth]{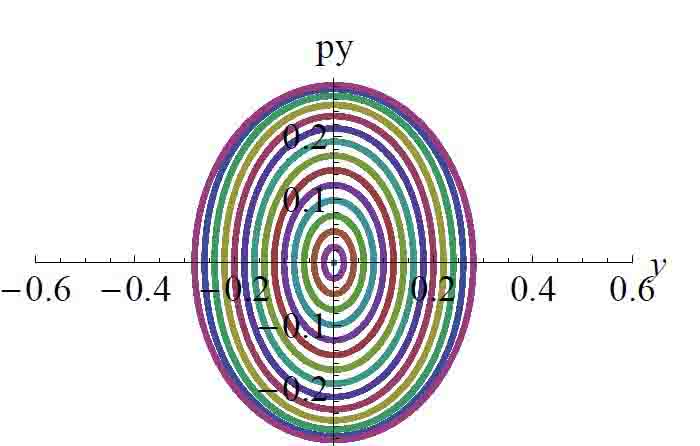}
\includegraphics[width=0.35\linewidth]{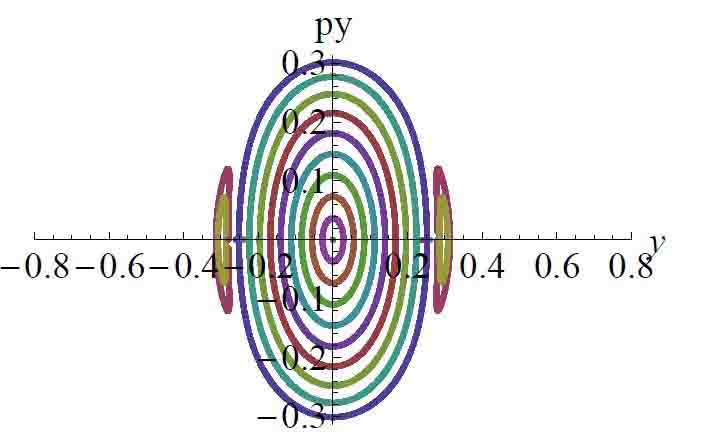}    \\
\includegraphics[width=0.35\linewidth]{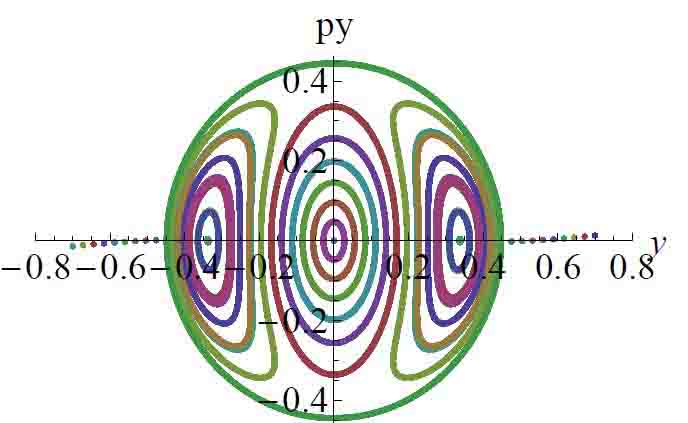}
\includegraphics[width=0.35\linewidth]{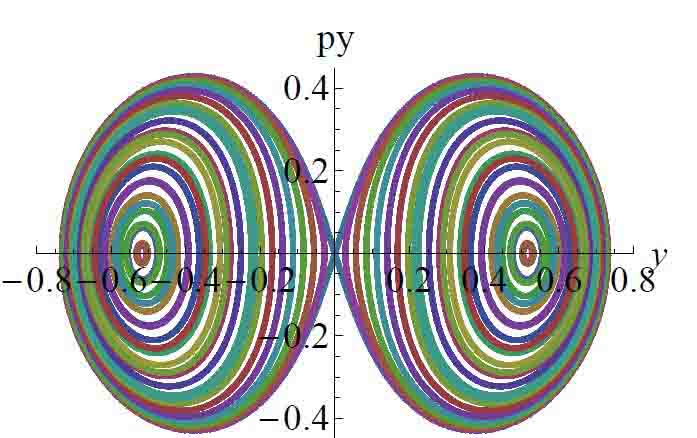}    \\
\caption{Poincar\'e sections of the center manifold associated to $L_1$ of the Sun--Vesta system on the plane $(y,p_{y})$
with $\beta=10^{-2}$ and $A=4.54776\cdot10^{-14}$;
different values of the energy are taken into account: $h=0.04$ (upper left panel), 0.05 (upper right), 0.1
(bottom left), 0.4 (bottom right).}\label{fig:Vestaconq}
\end{figure}

A different situation occurs when the radiation pressure and the oblateness are switched on, as shown in
Figure~\ref{fig:Vestaconq} which reports the Poincar\'e sections for the case with $\beta=10^{-2}$
and $A=4.54776\cdot10^{-14}$. Indeed, we have noticed that the oblateness has a small effect,
while the parameter $\beta$ plays a major r\^ole in shaping the dynamics. In fact, a comparison
between Figure~\ref{fig:Vestasenzaq} and Figure~\ref{fig:Vestaconq} shows that the sequence of
bifurcations is completely different.

In Figure~\ref{fig:Vestaconq} we observe a regular
behavior for $h=0.04$, while already at $h=0.05$ the dynamics experiences a first bifurcation
with the appearance of halo orbits and simultaneous loss of stability of
the planar Lyapunov orbit.
A second bifurcation of unstable inclined orbits takes place at about $h=0.1$; at this stage, the
planar Lyapunov orbit regains stability
(compare also with \cite{GM}), as shown in Figure~\ref{fig:Vestaconq} where the planar Lyapunov orbit is given
by the outermost curve.
For increasing values of the energy, the unstable families (which are located on the
vertical axis of the plot for $h=0.4$) collapse on the vertical Lyapunov orbit at the
center of the axes (see the bottom right panel of Figure~\ref{fig:Vestaconq}). This corresponds to
the third bifurcation, which is shown at about $h=0.4$.


\vskip.1in

\vskip.1in

\begin{figure}
\centering%
\includegraphics[width=0.28\linewidth]{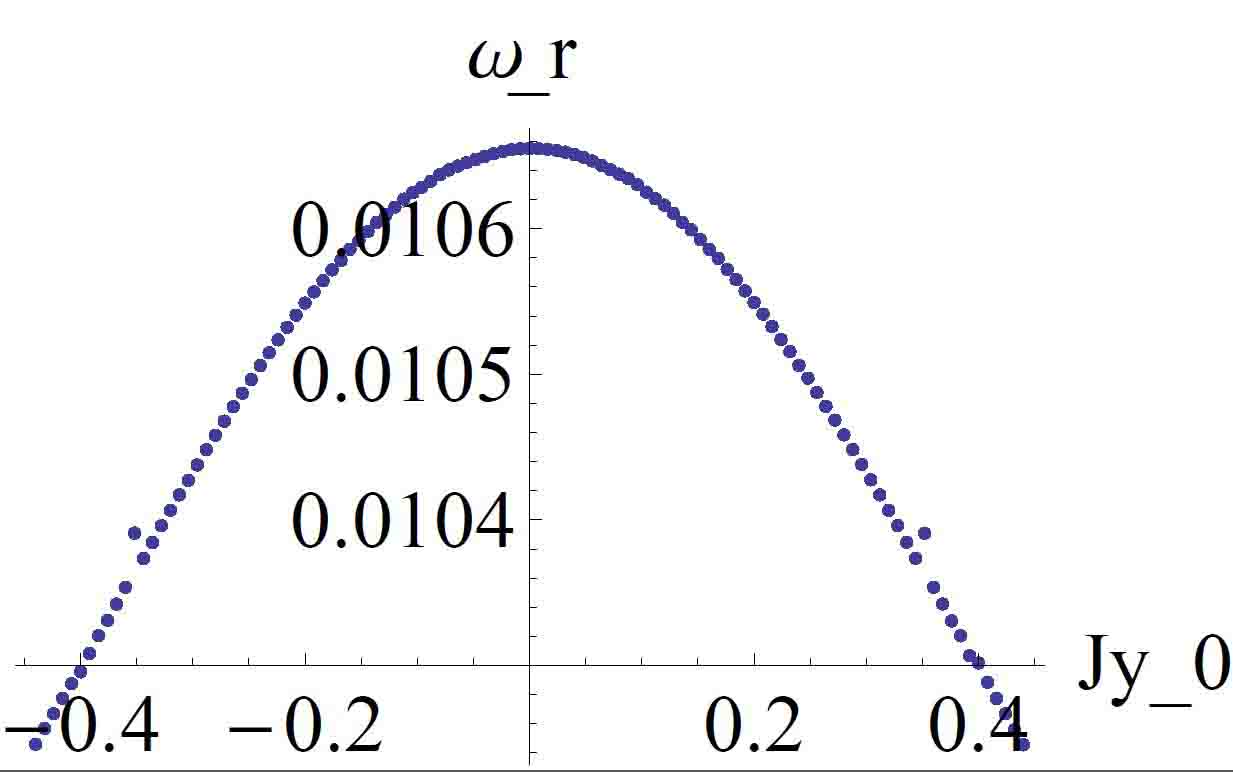}
\hglue0.2cm
\includegraphics[width=0.28\linewidth]{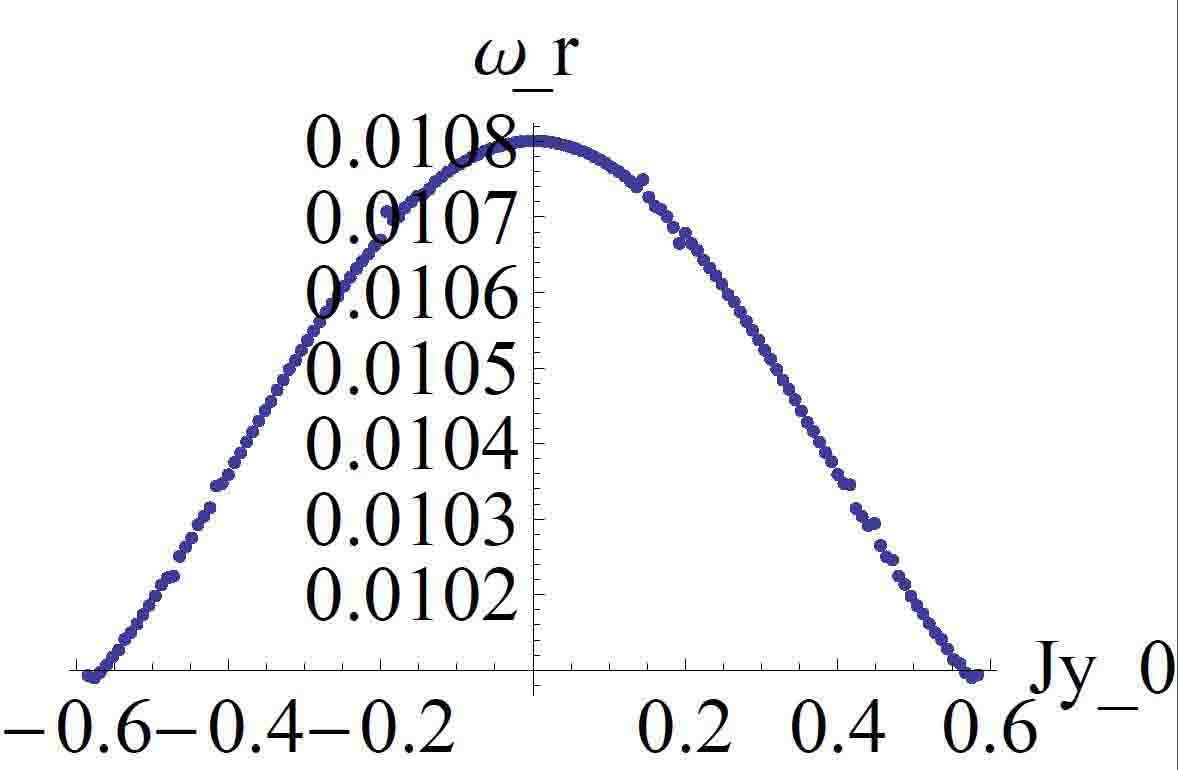}
\hglue0.2cm
\includegraphics[width=0.28\linewidth]{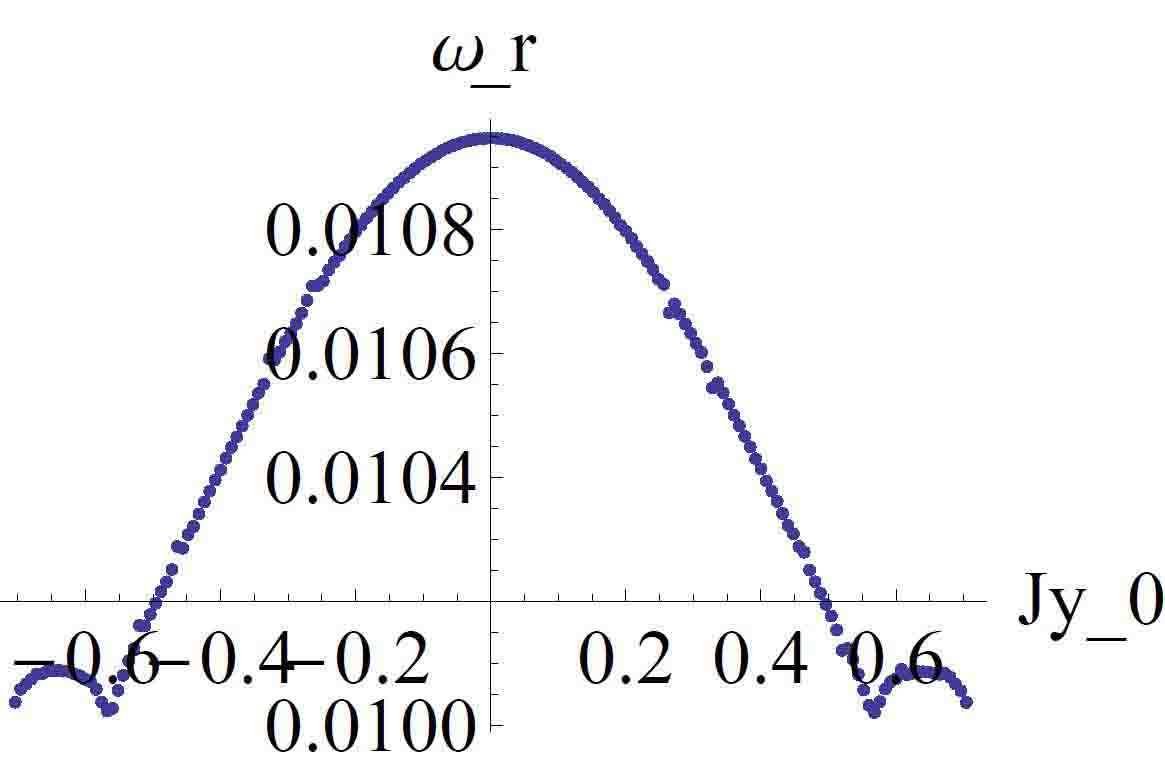}
\caption{Frequency analysis in the plane $(J_y^0,\omega_r)$ for the Sun--Vesta case with $\beta=0$, $A=0$;
different values of the energy are taken into account: $h=0.2$ (left panel), 0.35 (middle), 0.5 (right).}\label{fig:freqan}
\end{figure}

\begin{figure}
\centering
\includegraphics[width=0.35\linewidth]{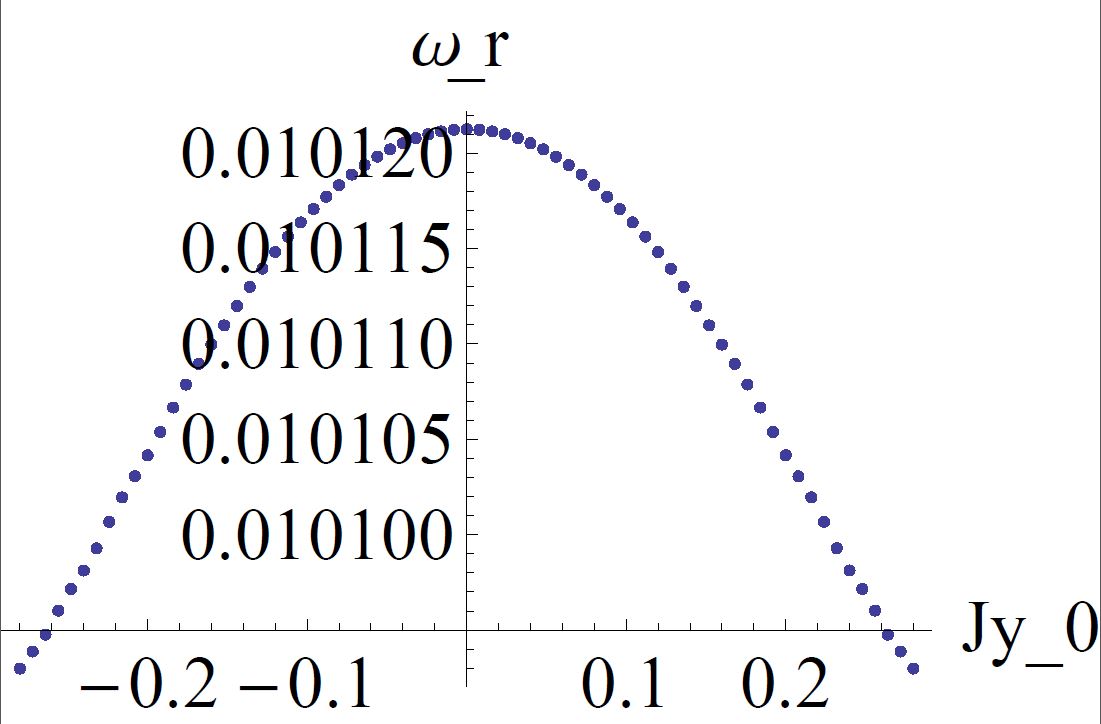}
\includegraphics[width=0.35\linewidth]{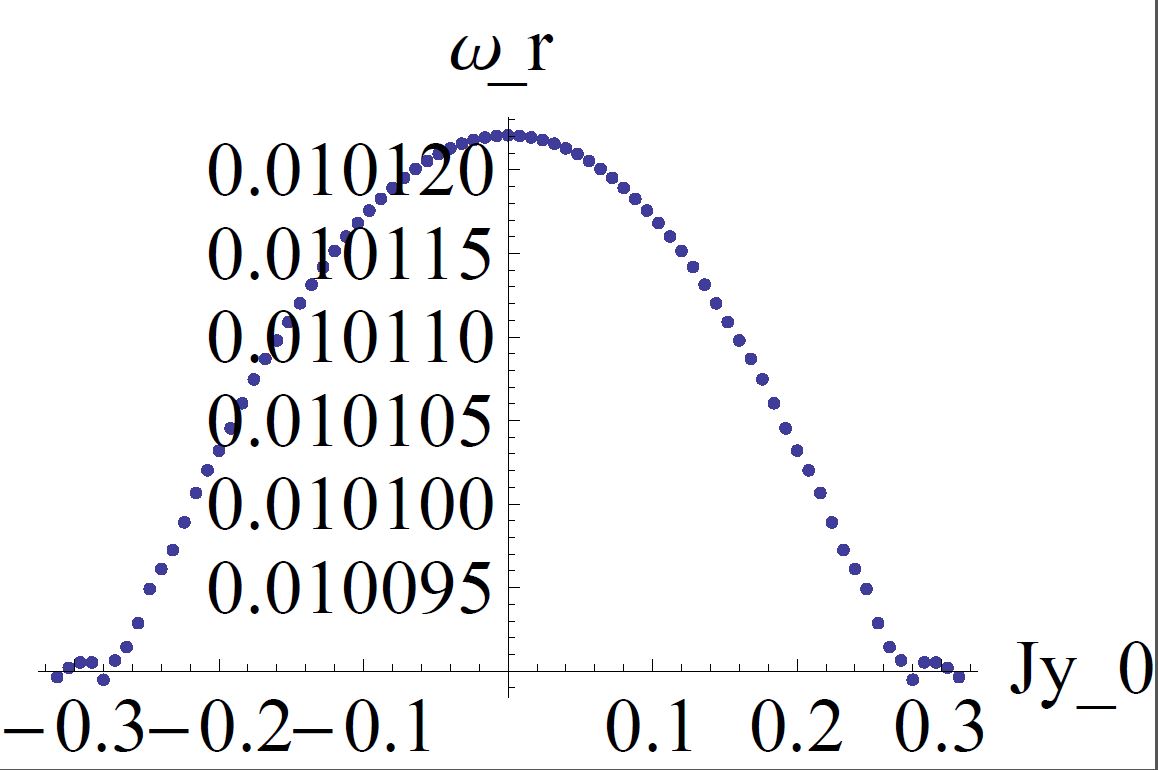}\\
\includegraphics[width=0.35\linewidth]{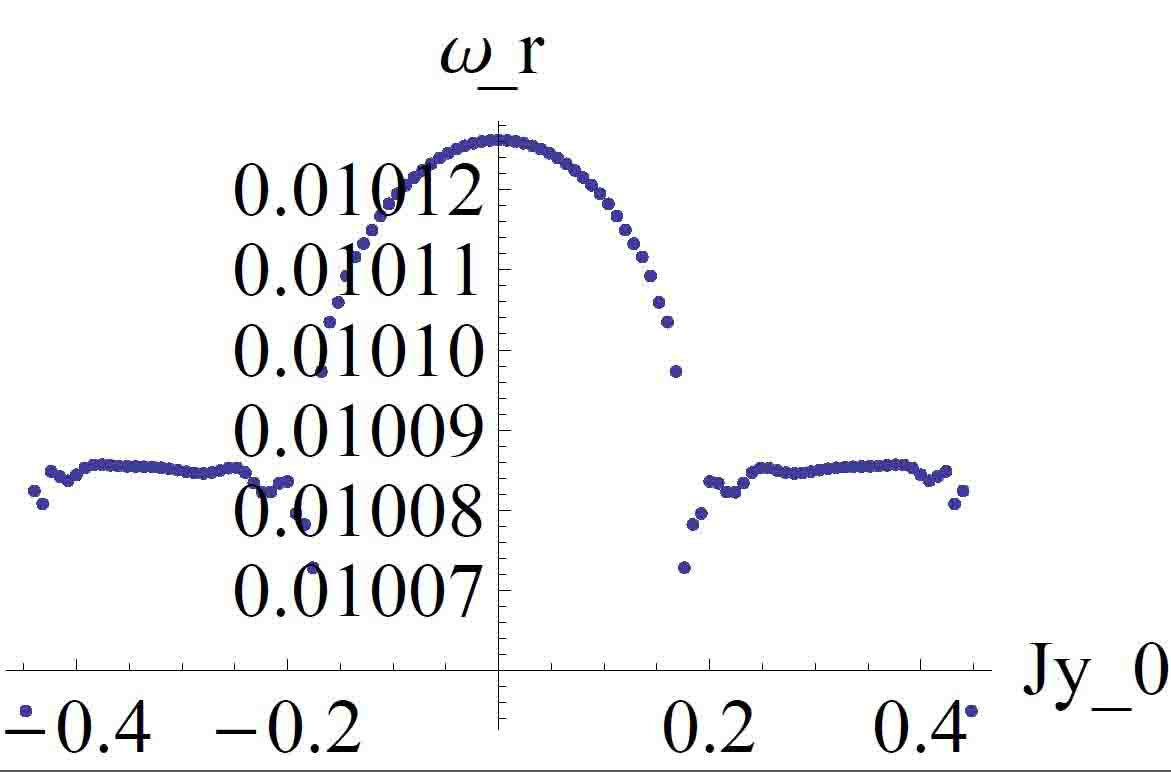}
\includegraphics[width=0.35\linewidth]{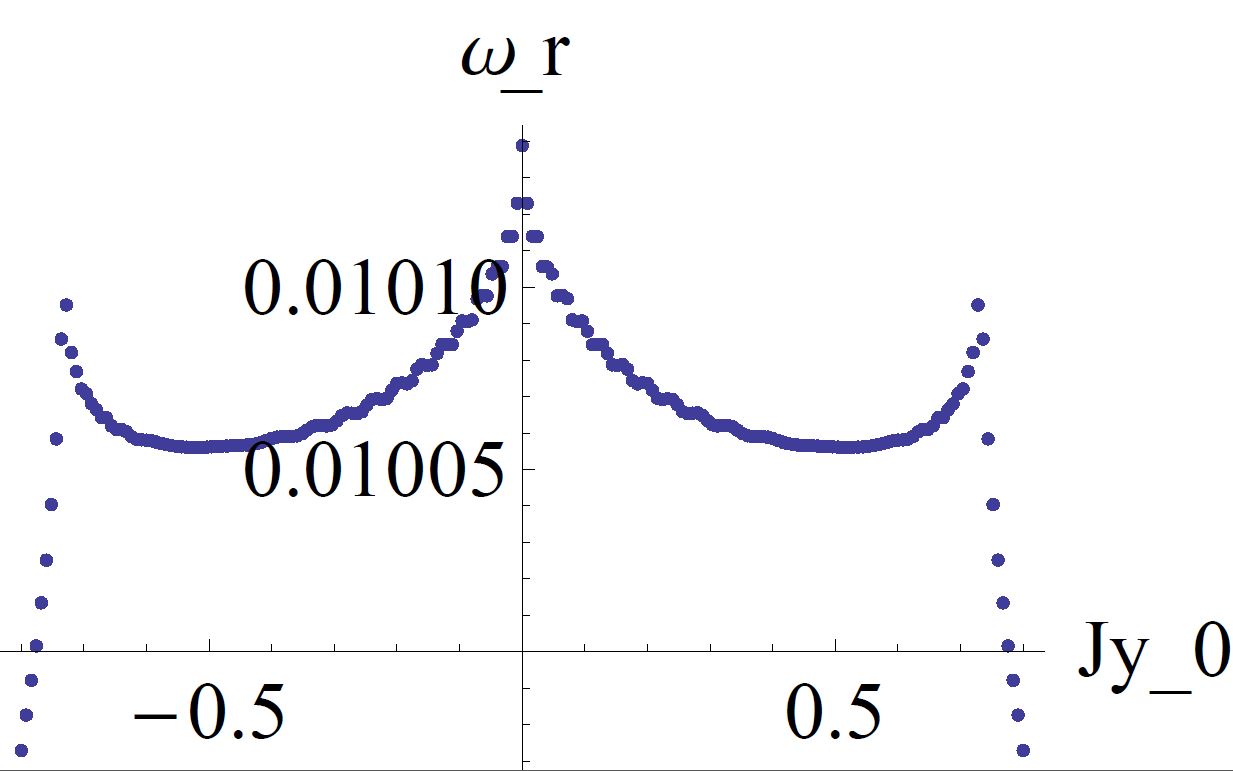}
\caption{Frequency analysis in the plane $(J_y^0,\omega_y)$ for the Sun--Vesta case with $\beta=10^{-2}$ and $A=4.54776\cdot10^{-14}$;
different values of the energy are taken into account: $h=0.04$ (upper left panel), 0.05 (upper right), 0.1
(bottom left), 0.4 (bottom right).}\label{fig:freqanq}
\end{figure}

The above results are confirmed by the study of the model through frequency analysis, as shown in
Figures~\ref{fig:freqan} and \ref{fig:freqanq}.

In particular, in Figure~\ref{fig:freqan} we provide the results of the frequency analysis for the Sun--Vesta case
without solar radiation pressure and no oblateness in the $(J_y^0,\omega_r)$ plane with $J_y^0$ the initial condition and
$\omega_r=|\omega_y/\omega_z|$. The first plot corresponds to $h=0.2$
and it shows a regular behavior as it was found in the first plot of Figure~\ref{fig:Vestasenzaq}. The occurrence
of small halo orbits in the frequency analysis investigation corresponds to the two tiny bumps at the outermost
sides of the plot for $h=0.35$ in Figure~\ref{fig:freqan}; these bumps increase in size for $h=0.5$, in full
agreement with the corresponding plot of Figure~\ref{fig:Vestasenzaq}.

\vskip.1in

The results of the investigation through frequency analysis in the case with solar radiation pressure and
oblateness are provided in Figure~\ref{fig:freqanq}, which corresponds to the Sun--Vesta case
with $\beta=10^{-2}$ and $A=4.54776\cdot10^{-14}$. Again we find full agreement with the Poincar\'e maps
provided in Figure~\ref{fig:Vestaconq}. Precisely, the upper left panel of Figure~\ref{fig:freqanq} shows
a regular behavior, while tiny bumps, corresponding to the bifurcation of the halo orbits for $h=0.05$,
are present in the upper right panel of Figure~\ref{fig:freqanq}. The three island regimes occurring
for $h=0.1$ correspond to the central bump and the two left and right wings of the bottom left panel.
Finally, for $h=0.4$ we observe a singular behavior on the vertical axis of the bottom right panel
of Figure~\ref{fig:freqanq}, which corresponds to the vertical Lyapunov orbit at the origin of
the coordinates in Figure~\ref{fig:Vestaconq}.

\vskip.1in

\begin{figure}
\centering
\includegraphics[width=0.35\linewidth]{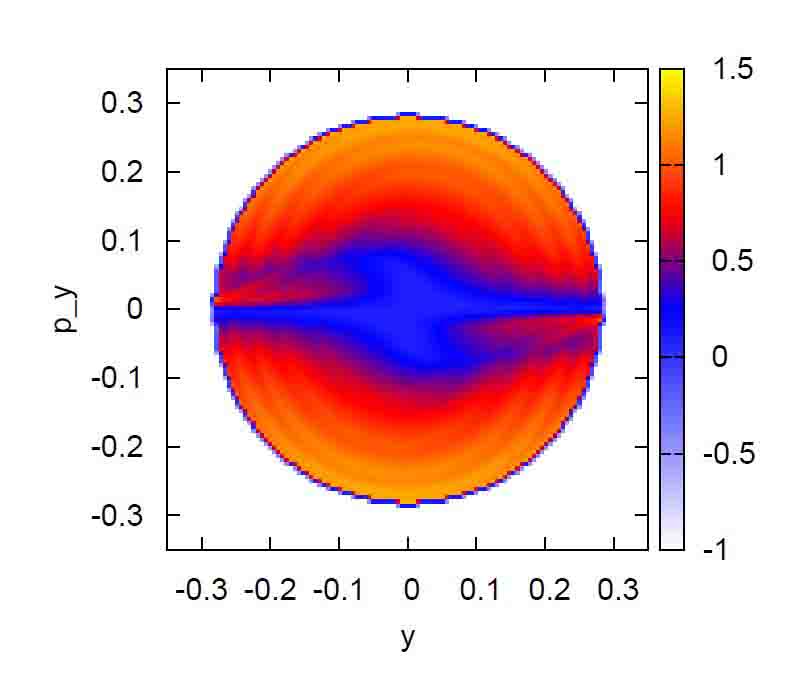}
\includegraphics[width=0.35\linewidth]{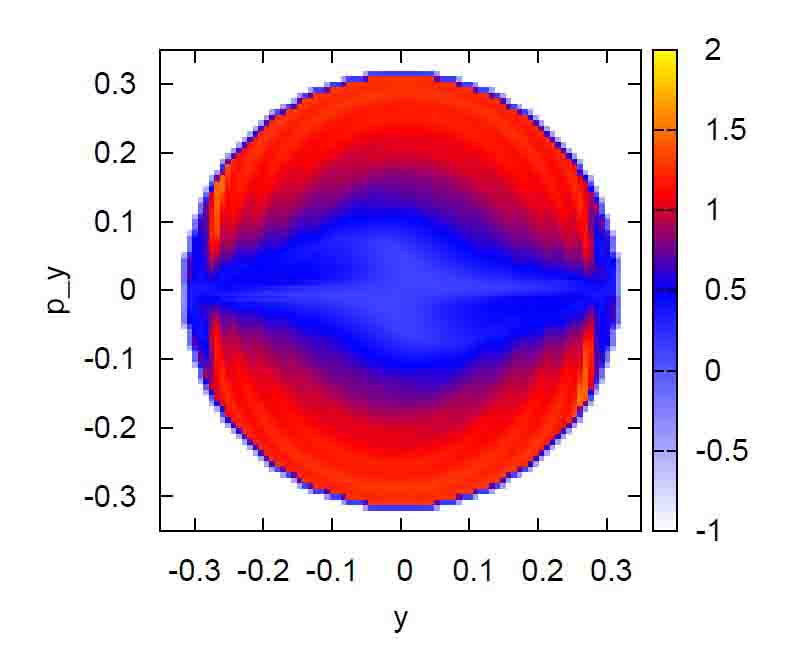}\\
\includegraphics[width=0.35\linewidth]{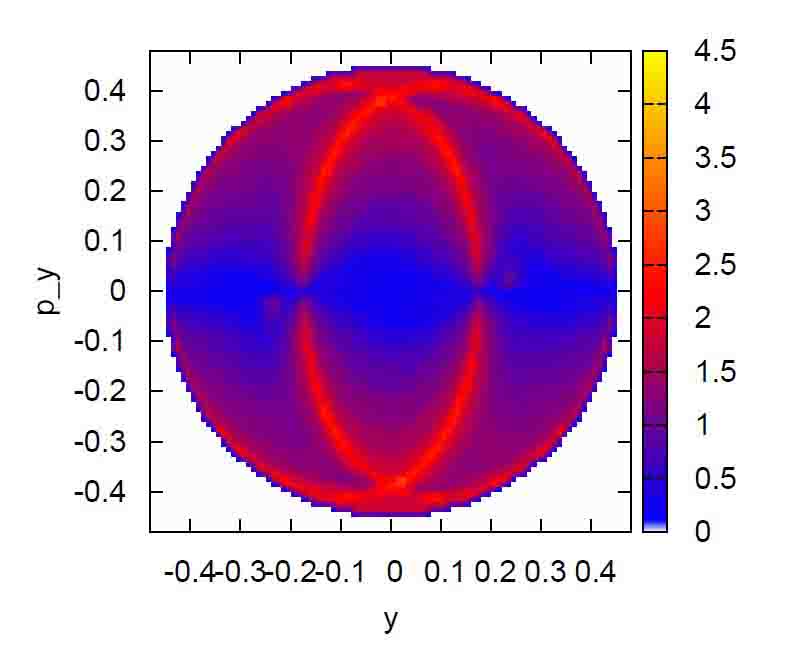}
\includegraphics[width=0.35\linewidth]{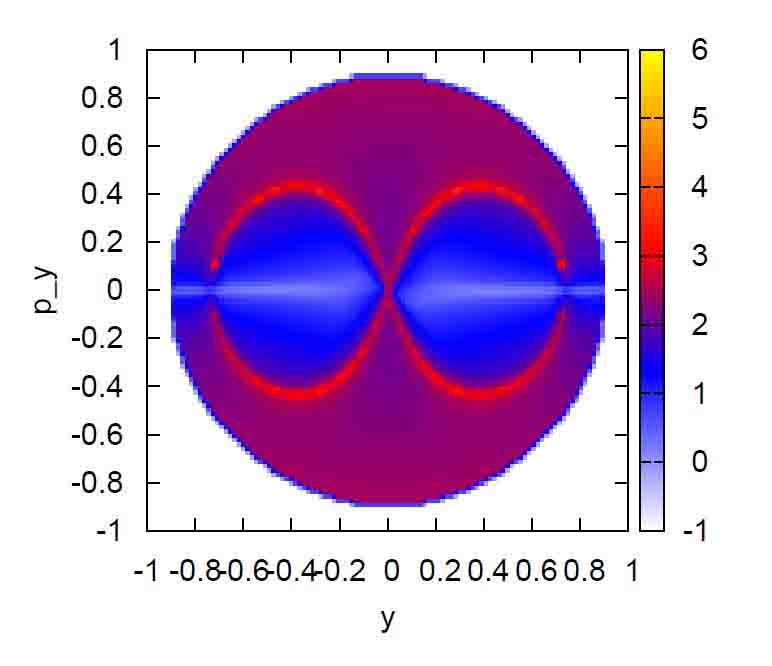}
\caption{Fast Lyapunov Indicators for the Sun--Vesta
case with $\beta=10^{-2}$ and $A=4.54776\cdot10^{-14}$; different
values of the energy are taken into account: $h=0.04$ (upper left
panel), 0.05 (upper right), 0.1 (bottom left), 0.4 (bottom
right).}\label{fig:FLIq}
\end{figure}

The computation of the FLIs provides additional information: beside the overall structure of the
phase space, it yields the regular or chaotic character of the different trajectories.
In particular, we can locate the separatrices, which were not easily determined within the
Poincar\'e maps (compare, e.g., Figure~\ref{fig:Vestaconq} bottom left and Figure~\ref{fig:FLIq}
bottom left).\\

\begin{figure}
\centering%
\includegraphics[width=0.28\linewidth]{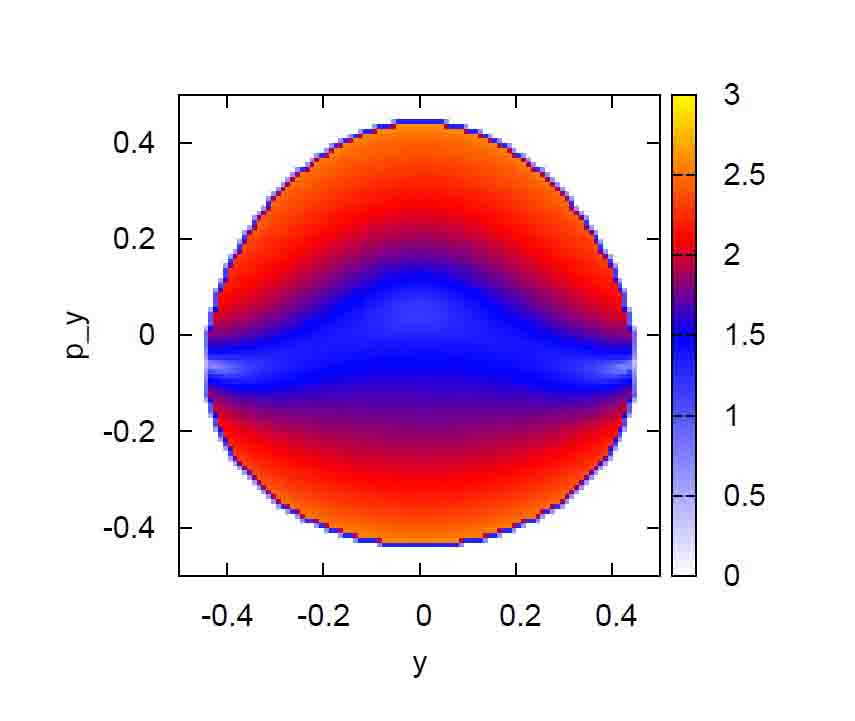}
\includegraphics[width=0.28\linewidth]{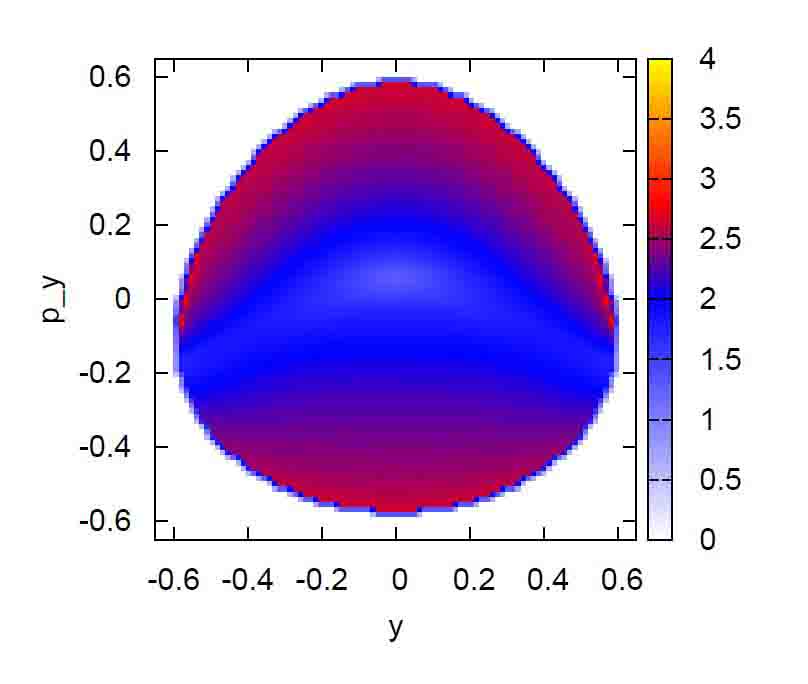}
\includegraphics[width=0.28\linewidth]{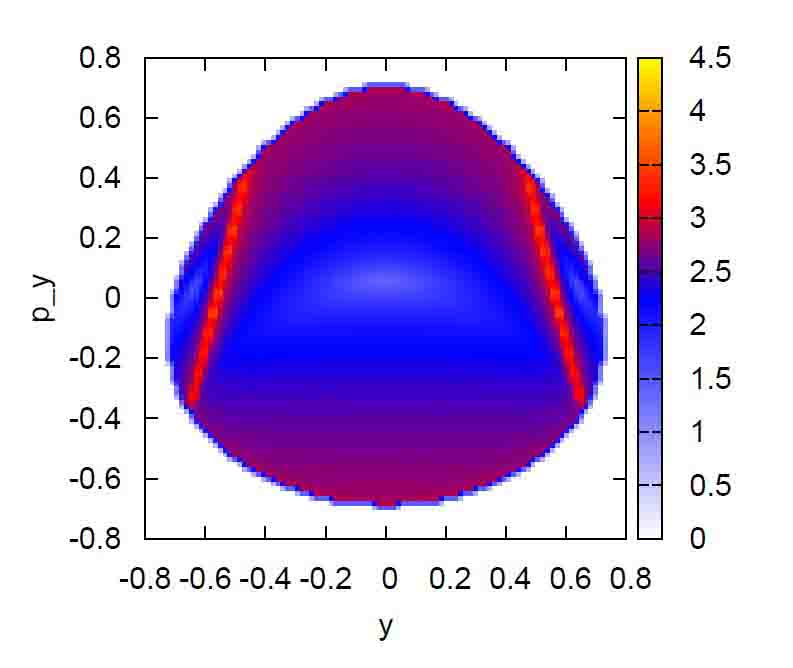}
\caption{Fast Lyapunov Indicators for the Sun--Vesta
case with $\beta=0$, $A=0$; different values of the energy are
taken into account: $h=0.2$ (left panel), 0.35 (middle),
0.5 (right).}\label{fig:FLI}
\end{figure}

In Figure~\ref{fig:FLI} we provide the results obtained computing the FLIs for
the Sun--Vesta case without solar radiation pressure and no oblateness; the results
must be compared with Figure~\ref{fig:Vestasenzaq} in order to distinguish the
different orbits on the Poincar\'e maps and to determine their stability on the
FLI plots. The color bar on the side of each plot gives the quantitative value of the FLI.

\vskip.1in

\vskip.1in

In Figure~\ref{fig:FLIq} we provide the results obtained computing the FLIs for the Sun--Vesta case
with $\beta=10^{-2}$ and $A=4.54776\cdot10^{-14}$. Again, Figure~\ref{fig:FLIq} must be compared with
Figure~\ref{fig:Vestaconq} in order to identify the various trajectories
on the Poincar\'e sections and to characterize their stability on the FLI plots.
It is remarkable how the FLI plots highlight also the separatrices as shown in particular in
the bottom panels of Figure~\ref{fig:FLIq}.

\vskip.1in

\begin{figure}
\centering
\includegraphics[width=0.35\linewidth]{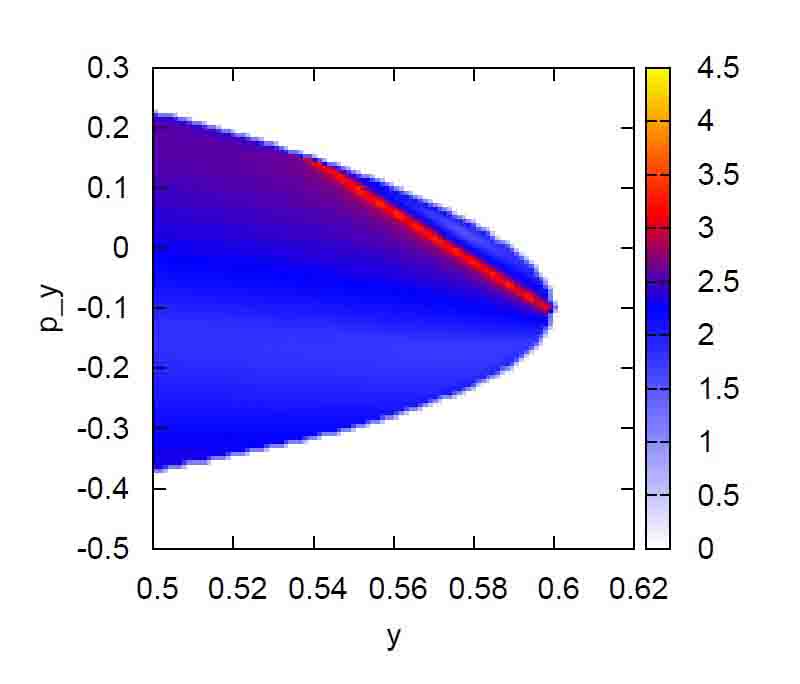}
\includegraphics[width=0.35\linewidth]{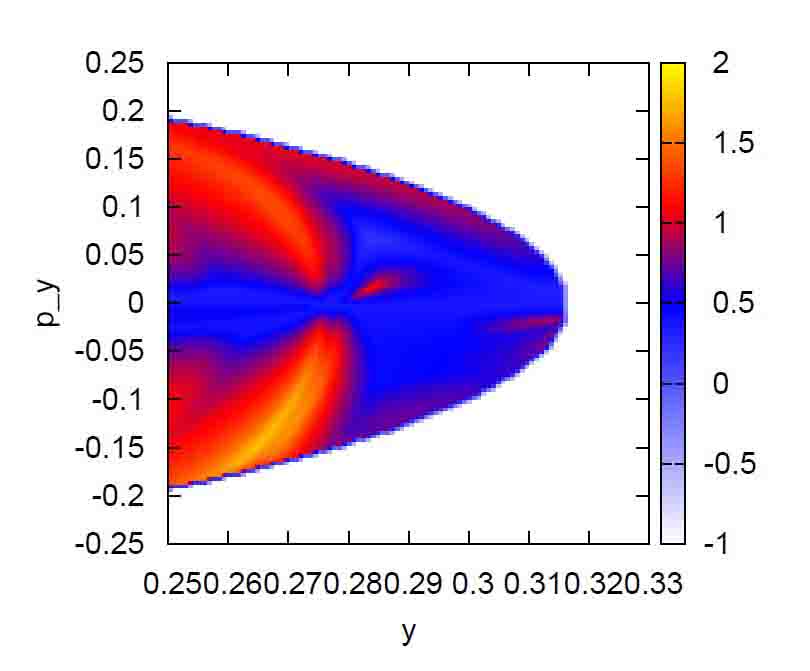}
\caption{A zoom on the cases with $\beta=0$, $h=0.35$ (left panel) and $\beta=10^{-2}$, $h=0.05$ (right).}\label{fig:FLIzoom}
\end{figure}

When the halo orbits are very tiny, a zoom becomes necessary, as shown in
Figure~\ref{fig:FLIzoom}, where we provide a magnification of the plots obtained in the cases
with $\beta=0$, $h=0.35$ (left panel) and $\beta=10^{-2}$, $h=0.05$ (right panel).

In the panels of Figure~\ref{fig:FLIq} we notice some lighter regions which are mainly along the horizontal
direction (compare with the bottom right panel); these zones do not have a real dynamical meaning,
but they are rather an artefact of the choice of the initial tangent vector, as the FLI strongly
depends on this choice. In Figure~\ref{fig:FLIq} the tangent vector has been fixed as $(v_{p_x},v_{p_y},v_x,v_y)=(1,0,0,0)$ and we observe that
lighter regions occur in the direction perpendicular to the chosen tangent vector. A reliable description of the dynamical character
of a specific region by means of the FLIs can only be obtained by comparing the plots produced using different
tangent vectors in orthogonal directions or, alternatively, by increasing very much the accuracy of the computations.
However, we remark that the analysis of Figure~\ref{fig:FLIq} suffices to distinguish the bifurcations as well as the main structures, like halo orbits or separatrices; further refinements go beyond the aims of the present work.

\section{Analytical versus numerical results}\label{sec:conclusion}
In this section we compare the results which are obtained implementing the analytical formulae \equ{inclined}-\equ{loop}
for the computation of the bifurcation thresholds with the numerical values obtained through the Poincar\'e maps or
the FLIs.

The analytical and numerical bifurcation values for all three case studies (Earth--Moon, Sun--barycenter, Sun--Vesta)
are listed in Table~\ref{tab:bif}. The analytical results require the computation of the normal form as well
as the computation of the thresholds at first order in the detuning as given by \equ{inclined}-\equ{loop}.
The qualitative results are obtained
looking at the bifurcations observed on the Poincar\'e sections as well as through the FLI maps
(the results are also validated through the application of the frequency analysis).

\vskip.1in

\begin{table}
\caption{Numerical (num) and analytical (anal) value of the energy at which the bifurcation to halo orbits has taken place; the results are given for the Earth--Moon, Sun--barycenter, Sun--Vesta systems.}\label{tab:bif}
\center
\tabcolsep=3mm
\renewcommand\arraystretch{1.2}
\begin{tabular}{|c|c|c|c|c|}
\hline
\hline
& $L_1$ num & $L_1$ anal & $L_2$ num & $L_2$ anal\\
\hline
\hline
Earth--Moon $\beta=0,A=0$ & $0.3026$ & $0.3069$ & $0.3776$ & $0.3636$\\
\hline
Earth--Moon $\beta=0,A=4.15559\cdot 10^{-9}$ & $0.3026$ & $0.3069$ & $0.3776$ & $0.3636$ \\
\hline
Sun--barycenter $\beta=0,A=0$ & $0.3333$ & $0.3356$ & $0.3377$  & $0.3391$ \\
\hline
Sun--barycenter $\beta=10^{-2},A=1.96782\cdot 10^{-12}$ & $0.2793$  & $0.2864$& $0.3759$  & $0.3755$\\
\hline
Sun--Vesta $\beta=0,A=0$ & $0.3341$  & $0.3373$& $0.3346$  & $0.3374$\\
\hline
Sun--Vesta $\beta=10^{-2},A=4.54776\cdot 10^{-14}$ & $0.0422$  & $0.0424$& $0.4451$  & $0.4434$\\
\hline
\end{tabular}
\end{table}

\vskip.1in

From Table~\ref{tab:bif} we can draw the following conclusions.

\begin{itemize}
    \item[$(i)$] The agreement between the analytical and numerical results is very satisfactory.
    The relative error between the theoretical and experimental values ranges between
    $10^{-2}$ and $10^{-3}$.
    \item[$(ii)$] A first order estimate is already enough to get a good approximation of
    the bifurcation thresholds; this estimate requires a very little computational effort with
    respect to the qualitative analysis based on the Poincar\'e maps or the FLIs. Obviously,
    better results can be obtained computing higher order normal forms, but at the expense
    of dealing with more complex formulae.
    \item[$(iii)$] Switching on the solar radiation pressure provokes drastic changes for small mass parameters. In particular,
in the Sun--Vesta case the first, second and third bifurcations take place at much lower values of the energy level,
such that the other bifurcations become feasible. From the physical point of view, the reason for such a peculiar behavior is
due to the balance between a smaller mass like that of an asteroid and the effect of the solar radiation pressure.
    \item[$(iv)$] The r\^{o}le of the oblateness is essentially negligible in all considered cases. This fact
    could have been inferred easily, but we believe worthwhile to derive a complete model, valid not only for the
    cases studied in the present paper, but also for general situations in which the small body could have
    a very irregular shape. Simple experiments show that in a Sun--asteroid sample, the oblateness becomes important only
    when the factor $A$ is as large as $10^{-6}-10^{-7}$. This parameter value does not apply to Vesta, but it
    might be of interest for other astronomical situations.
   \end{itemize}

\section*{Acknowledgments}

We deeply thank Gerard Gomez for very useful suggestions and Danilo Stella for
interesting discussions.
M.C. and A.C. were partially supported by the European Grant MC-ITN Astronet-II,
G.P. was partially supported by the European Grant MC-ITN Stardust. M.C., A.C., G.P. acknowledge GNFM/INdAM.

\appendix

\section[]{Derivation of the mean motion around an oblate primary}\label{app:meanmotion}
For completeness we derive here the formula for the mean motion of a massless body $P$ under the gravitational effect of an oblate primary.
Let us denote by $M_P$ and $A$ the mass and the oblateness coefficient of the primary, $r$ represents the distance of $P$
from the primary and $\mathcal{H}$ is the angular momentum constant. Then, the effective potential (see, e.g., \cite{alebook}) can be written as
$$
V_{eff}(r)={\mathcal{H}^2\over{2r^2}}-{M_P\over r}-{{M_P\, A}\over {2r^3}}\ ,
$$
whose derivative is
$$
V_{eff}'(r)=-{\mathcal{H}^2\over{r^3}}+{M_P\over r^2}+{{3M_P\,A}\over {2r^4}}\ .
$$
The solutions of $V_{eff}'(r)=0$ are given by
\beq{rapp}
r={{\mathcal{H}^2\pm\sqrt{\mathcal{H}^4-6M_P^2\,A}}\over {2M_P}}\ .
\eeq
Taking into account that $A$ is small, we can approximate the non--trivial solution of \equ{rapp} by
\beq{rhapp}
r={\mathcal{H}^2\over M_P}\ (1-{3\over 2}{{M_P^2\,A}\over \mathcal{H}^2})\ .
\eeq
Assume that the orbit of $P$ is circular, say $r=a$, we have that $\mathcal{H}^2=n^2a^4$ which, together with
\equ{rhapp} and the normalization of the units of measure such that $a=1$, $M_P=1$, provides:
$$
n^2=1+{3\over 2}A\ .
$$

\section[]{Reduction of the quadratic part}\label{app:diag}
In this section we provide the details of the reduction of the quadratic part as performed in
Section \ref{sec:reduction} in order to obtain the Hamiltonian \equ{eq5} (equivalently \equ{eq6}).
The procedure is very similar to that explained in \cite{JM}, to which we refer for a complete discussion; however, for self-consistency,
we provide here some details containing the necessary amendments to encompass the oblate case with solar radiation pressure.\\

We start by computing the eigenvalues of $M$ in \equ{Mdef}; we denote by $I_n$ the $n$--dimensional identity matrix.
We notice that by defining $M_{\lambda}\equiv M-\lambda I_{4}$, we can write
$
M_{\lambda}=
\left(%
\begin{array}{cc}
A_{\lambda}& I_{2}\\
B & A_{\lambda}\\
 \end{array}%
\right)
$
with $A_{\lambda}=
\left(%
\begin{array}{cc}
-\lambda & n\\
-n & -\lambda\\
 \end{array}%
\right)
$ and
$
B=
\left(%
\begin{array}{cc}
-2a & 0\\
0 & -b\\
 \end{array}%
\right)$.
Then, the kernel of $M_{\lambda}$ is given by the solution of $M_{\lambda} w=0$
with
$w=
\left(
\begin{array}{c}
w_{1}\\
w_{2}\\
\end{array}
\right)$ and $w_{1}$,$w_{2}\in\bf{R}^{2}$.
Simple computations show that the eigenvector of $M$ is given by
$$
(2n\lambda ,\  \lambda^{2}+2a-n^{2}, \ n \lambda^{2}-2an+n^{3},\ \lambda^{3}+(2a+n^{2})\lambda)^\top,
$$
where $\lambda$ is an eigenvalue of $M$ (the superscript $\top$ denotes the transposed).\\

Let us consider the eigenvectors associated to $\omega_{1}=\sqrt{-\eta_{1}}$; from \equ{plam} we have $p(\lambda)=0$, so that
$\omega_{1}$ satisfies the equation
$$
\omega_{1}^{4}-(2n^{2}+2a+b)\omega_{1}^{2}+(n^{4}-2an^{2}-bn^{2}+2ab)=0\ .
$$
Using $i\omega_{1}=\lambda_{1}$, the eigenvector $u_{\omega_{1}}+iv_{\omega_{1}}$ associated to $\omega_{1}$ is given by
\beqano
u_{\omega_{1}}+iv_{\omega_{1}}&\equiv&(2ni\omega_{1},-\omega_{1}^{2}+2a-n^{2},\nonumber\\
&-&n\omega_{1}^{2}-2an+n^{3}, -i\omega_{1}^{3}+(2a+n^{2})i\omega_{1})^\top\ ,
\eeqano
while for $\pm\lambda_{1}=\sqrt{\eta_{2}}$ we obtain:
\beqano
u_{+\lambda_{1}}&=&(2n\lambda_{1},\lambda_{1}^{2}+2a-n^{2},\nonumber\\
&&n\lambda_{1}^{2}-2an+n^{3},\lambda_{1}^{3}+(2a+n^{2})\lambda_{1})^\top\nonumber\\
v_{-\lambda_{1}}&=&(-2n\lambda_{1},\lambda_{1}^{2}+2a-n^{2},\nonumber\\
&&n\lambda_{1}^{2}-2an+n^{3},-\lambda_{1}^{3}-(2a+n^{2})\lambda_{1})^\top\ .
\eeqano
Let $C=(u_{+\lambda_{1}},\ u_{\omega_{1}},\ v_{-\lambda_{1}},\  v_{\omega_{1}})$; we have
$$
C^\top JC=
\left(%
\begin{array}{cc}
0 & D\\
-D & 0\\
 \end{array}%
\right)
\ ,\quad
\mbox{where}\ \ \
D=
\left(%
\begin{array}{cc}
d_{\lambda_{1}} & 0\\
0 & d_{\omega_{1}}\\
 \end{array}%
\right)
$$
with $d_{\lambda_{1}}$ and $d_{\omega_{1}}$ given by
\beqano
d_{\lambda_{1}}&=&-2\lambda_{1}((-4n^{2}+2a-b)\lambda_{1}^{2}-4n^{4}+ bn^{2} + 6an^{2}-2ab+4a^{2})\ ,\nonumber\\
d_{\omega_{1}}&=&-\omega_{1}((-4n^{2}+2a-b)\omega_{1}^{2}+4n^{4}-bn^{2}-6an^{2}+2ab-4a^{2})\ .
\eeqano
In order to obtain a symplectic change of variables, we re--scale by $s_{1}=\sqrt{d_{\lambda_{1}}}$,
$s_{2}=\sqrt{d_{\omega_{1}}}$ and we require that $d_{\lambda_{1}}>0$, $d_{\omega_{1}}>0$.
Finally, we re--scale $(z,p_z)$ by $(\frac{1}{\sqrt{\omega_{2}}},\sqrt{\omega_{2}})$. The final symplectic
change of variables is given by the matrix whose columns are $u_{+\lambda_{1}}/s_1$,
$u_{\omega_{1}}/s_2$, ${v_{-\lambda_{1}}/s_1}$, $v_{\omega_{1}}/s_2$.

\section[]{Center manifold reduction}\label{sec:Lie}
Let $H$ be a Hamiltonian function with 3 degrees of freedom, admitting an equilibrium point of type
\sl saddle$\times$center$\times$center. \rm Let $\pm\lambda$, $\pm i \omega_{1}$, $\pm i\omega_{2}$ be the eigenvalues
of the linearized system. Expanding the Hamiltonian around the equilibrium point in complex variables,
we obtain a simpler Hamiltonian function of the form
$$
 H(q,p)= H_{2}(q,p)+ \sum_{n\geq 3}H_{n}(q,p)\ ,
$$
where
$$
H_{2}(q,p)=\lambda q_{1}p_{1}+i\omega_{1}q_{2}p_{2}+i\omega_{2}q_{3}p_{3}\ ,
$$
while $H_{n}(q,p)$ are homogeneous polynomials of degree $n$ in the variables $(q,p)$.
To decouple the hyperbolic direction from the elliptic one, we need to kill the monomials whose
exponent in $p_1$ is different from that in $q_1$. This procedure, which makes use of Lie series, allows one to get a first integral
with level surface given by the center manifold. We sketch below the procedure to find the canonical transformation,
referring to \cite{JM} for full details. \\

Given a Hamiltonian $H$ and a generating function $G$, we denote by $\hat{H}$ the function\footnote{Curly brackets denote, as usual, the Poisson brackets
(\cite{Goldstein}).}
\beq{eq9}
\hat{H} \equiv H + \left\lbrace H, G\right\rbrace  + \frac{1}{2!}\left\lbrace \left\lbrace H, G\right\rbrace , G\right\rbrace  + \frac{1}{3!}\left\lbrace \left\lbrace \left\lbrace H, G\right\rbrace , G\right\rbrace , G\right\rbrace  +...
\eeq
If $G$ has degree 3, say $G=G_3$, comparing same orders in \equ{eq9} provides
\beqano
\hat{H_{2}}&=& H_{2}\ ,\nonumber\\
\hat{H_{3}}&=& H_{3}+ \left\lbrace H_{2}, G_{3}\right\rbrace\ ,\nonumber\\
\hat{H_{4}}&=& H_{4}+  \left\lbrace H_{2}, G_{4}\right\rbrace\ +\left\lbrace H_{3}, G_{3}\right\rbrace +
\frac{1}{2!}\left\lbrace \left\lbrace H_{2}, G_{3}\right\rbrace , G_{3}\right\rbrace\ ,\quad ...
\eeqano
Next we look for $G_{3}$ such that $\hat{H_{3}}$ is in normal form. Expanding $H_{2}$, $H_{3}$, $G_{3}$ as
\beqano
H_{2}(q,p)&=&\sum_{j=1}^{3}\eta_{j}q_{j}p_{j}\ ,\nonumber\\
H_{3}(q,p)&=&\sum_{|k_{q}|+|k_{p}|=3}h_{k_{q},k_{p}}q^{k_{q}}p^{k_{p}}\ ,\nonumber\\
G_{3}(q,p)&=&\sum_{|k_{q}|+|k_{p}|=3}g_{k_{q},k_{p}}q^{k_{q}}p^{k_{p}}
\eeqano
for some coefficients $h_{k_{q},k_{p}}$, $g_{k_{q},k_{p}}$, $\eta_{j}$ and denoting by
$\eta=(\lambda, i\omega_{1},i\omega_{2})$, we have
$$
G_{3}(q,p)=\sum_{(k_{q}, k_{p}) \in \mathcal{S}_{3}}\frac{-h_{k_{q},k_{p}}}{<k_{p}-k_{q},\eta>}q^{k_{q}}p^{k_{p}}\ ,
$$
where $\mathcal{S}_{3}$ is the set of indexes $(k_{p}, k_{q})$, such that
$|k_{p}| + |k_{q}| = 3$ and with the first component of $k_{p}$ different from the first component of $k_{q}$.
We proceed iteratively to higher orders up to a given order, say $N$, so that we obtain:
$$
\hat{H}(q,p)=H^{(N)}(q_{1}p_{1}, q_{2}, q_{3}, p_{2}, p_{3})+ R^{(N+1)}(q_{1}, q_{2}, q_{3}, p_{1}, p_{2}, p_{3})\ ,
$$
where $H^{(N)}$ is a polynomial of degree $N$ and $R^{(N+1)}$ is a reminder of order $N+1$. Neglecting the reminder
and setting $q_{1}p_{1}=0$, we eliminate the hyperbolic component of  $H^{(N)}$ and we obtain a 2 degrees of freedom Hamiltonian
of the desired form. As remarked in \cite{JM}, there are no small divisors in the above procedure, since
$|<k_{q}-k_{p},\nu>| \geq |\lambda_1|$ for any $(k_{p}, k_{q})\in \mathcal{S}_j$, $j\geq 3$.\\

Going back to real variables $(y,z,p_{y},p_{z})$, we obtain a Hamiltonian of the form
$$
\tilde H(y,z,p_{y},p_{z})=\sum_{k_1,k_2,k_3,k_4\in{\bf Z}} h_{k_1,k_2,k_3,k_4}\ y^{k_1} z^{k_2} p_y^{k_3} p_z^{k_4}\ .
$$
The first few non--zero terms $h_{k_1,k_2,k_3,k_4}$ of the Hamiltonian restricted to the center manifold are provided
in Table~\ref{tab:hHamSV}.

\vskip.1in

\begin{table}
\caption{Coefficients up to degree 4 of the Hamiltonian restricted to the center manifold for the Sun--Vesta system
with $\beta=10^{-2}$. The exponents
$(k_{1}$, $k_{2}$, $k_{3}$, $k_{4})$ refer to the variables $(y,z,p_{y},p_{z})$.}\label{tab:hHamSV}

\center
\tabcolsep=3mm
\begin{tabular}{|c|c|c|c|c|}
  \hline
  $k_{1}$ &  $k_{2}$ & $k_{3}$ &  $k_{4}$ & $h_{k}$ \\
  \hline
  \hline
  2 & 0& 0& 0&0.501797549378742 \\
  \hline
  0 & 2 & 0 & 0 & 0.500906031584819 \\
  \hline
  0 & 0 &2 & 0&0.501797549378742 \\
  \hline
  0 & 0 &0 & 2&0.500906031584819\\
  \hline
  2 & 0 &1 & 0&0.0014920550494420622 \\
  \hline
  0 & 0 &3 & 0& -0.000248666270046148 \\
  \hline
  0 & 2 &1 & 0&0.0003790464810461912 \\
  \hline
  4 & 0 &0 & 0&  -0.02099512477285749\\
  \hline
  2 & 2 &0 & 0&-0.010667382077111268\\
  \hline
  0 & 4 &0 & 0&   -0.001354993618855492\\
  \hline
  2 & 0 &2 & 0& 0.04199066782897781 \\
  \hline
  0 & 2 &2 & 0&0.010667253491139606\\
  \hline
  0 & 0 &4 & 0&  -0.003498979471471415\\
  \hline
  1 & 1 &1 & 1&$2.81508155360122\cdot 10^{-7}$ \\
  \hline
  2 & 0 &0 & 2&$1.8824017340799554\cdot 10^{-7}$ \\
  \hline
  0 & 2 &0 & 2&  $4.7821141283422436\cdot 10^{-8}$ \\
  \hline
  0 & 0 &2 & 2&$-9.411646402154861\cdot 10^{-8}$ \\
  \hline
 \end{tabular}
\end{table}

\end{document}